\documentclass[11pt]{article}
\usepackage{amsmath, amssymb, amsfonts}
\usepackage{epsfig}

\def\be{\begin{equation}}
\def\ee{\end{equation}}
\def\bea{\begin{eqnarray}}
\def\eea{\end{eqnarray}}
\def\bes{\begin{eqnarray*}}
\def\ees{\end{eqnarray*}}

\def\nn{\nonumber}
\def\<{\langle}
\def\>{\rangle}
\def\lb{\label}
\def\bs{\setminus}

\def\R{{\mathbb{R}}}

\def\Z{{\mathbb{Z}}}
\def\N{{\mathbb{N}}}

\def\Q{{\mathbb{Q}}}

\def\RP{{\mathbb{R}P^{2n+1}}}

\def\aa{{\alpha}}
\def\bb{{\beta}}
\def\ga{{\gamma}}

\def\om{{\omega}}

\def\Lm{{\Lambda}}

\def\LS{{\cal LS}}

\def\ev{{\it ev}}
\def\id{{\it id}}

\def\ker{{\rm ker}}

\def\hv{{\rm hv}}
\def\diag{{\rm diag}}

\def\rank{{\rm rank}}

\def\mod{{\rm mod}}
\def\CG{{\rm CG}}
\def\Span{{\rm Span}}

\def\ol#1{\overline{#1}}  %overline in math mode
\def\td#1{\tilde{#1}}

\def\mapright#1{\smash{\mathop{\longrightarrow}\limits^{#1}}}

\def\mapdown#1{\Big\downarrow\rlap{$\vcenter{\hbox{$\scriptstyle#1$}}$}}

\def\mapright#1{\smash{\mathop{\longrightarrow}\limits^{#1}}}

\def\mapdown#1{\Big\downarrow\rlap{$\vcenter{\hbox{$\scriptstyle#1$}}$}}

\def\mapse#1{\searrow\rlap{$\vcenter{\hbox{$\scriptstyle#1$}}$}}
\def\mapsw#1{\swarrow\rlap{$\vcenter{\hbox{$\scriptstyle#1$}}$}}

\title{Topological structure of non-contractible loop space and\\
closed geodesics on real projective spaces with odd dimensions}
\author{Yuming Xiao $^{1}$\thanks{Partially supported by the Funds for Young Teachers of Sichuan University.
e-mail: yumingxiao@scu.edu.cn} \quad and \quad
Yiming Long$^{2}$\thanks{Partially supported by NSFC Grant 11131004, MCME, LPMC of MOE of China,
Nankai University and BCMIIS of Capital Normal University. e-mail: longym@nankai.edu.cn}\\\\
$^{1}$ School of Mathematics, Sichuan University,\\ Chengdu 610064, People's Republic of China \\
$^{2}$ Chern Institute of Mathematics and LPMC, Nankai University,\\ Tianjin 300071,  People's Republic of China\\
}
\begin{document}
\newtheorem{definition}{Definition}[section]
\newtheorem{theorem}{Theorem}[section]
\newtheorem{lemma}{Lemma}[section]
\newtheorem{corollary}{Corollary}[section]
\newtheorem{example}{Example}[section]
\newtheorem{property}{Property}[section]
\newtheorem{proposition}{Proposition}[section]
\newtheorem{remark}{Remark}[section]

\newcommand{\qed}{\nolinebreak\hfill\rule{2mm}{2mm}
\par\medbreak}
\newcommand{\Proof}{\par\medbreak\it Proof: \rm}
\newcommand{\rem}{\par\medbreak\it Remark: \rm}
\newcommand{\defi}{\par\medbreak\it Definition : \rm}
\renewcommand{\thefootnote}{\arabic{footnote}}
\maketitle

\begin{abstract}
{\it In this paper, we use Chas-Sullivan theory on loop homology and Leray-Serre spectral sequence to
investigate the topological structure of the non-contractible component of the free loop space on the
real projective spaces with odd dimensions. Then we apply the result to get the resonance identity of
non-contractible homologically visible prime closed geodesics on such spaces provided the total number
of such distinct closed geodesics is finite. }
\end{abstract}



\renewcommand{\theequation}{\thesection.\arabic{equation}}
\renewcommand{\thefigure}{\thesection.\arabic{figure}}

\setcounter{figure}{0}
\setcounter{equation}{0}
\section{Introduction and the main results}%Section 1

In this paper, we study the topological structure of the non-contractible component of the free loop space
on the real projective spaces with odd dimensions $\mathbb{R}P^{2n+1}$, which are the typically oriented and
non-simply connected manifolds with the fundamental group $\pi_{1}(\mathbb{R}P^{2n+1}) = \mathbb{Z}_{2}$.
Then we apply the result to get the resonance identity of non-contractible homologically visible prime closed
geodesics on such spaces when the total number of such distinct closed geodesics is finite.

Let $(M,F)$ be a Finsler manifold with finite dimension. As usual we choose an auxiliary Riemannian
metric on $M$ (cf. \cite{Shen2001}). It endows the free loop space $\Lm M$ on $M$ defined by
$$  \Lambda M=\left\{c: S^{1}\to M\mid c\ {\rm is\ absolutely\ continuous\ and}\
                        \int_{0}^{1}F(\dot{c},\dot{c})dt<+\infty\right\},  $$
with a natural structure of Riemannian Hilbert manifold on which the group $S^1=\R/\Z$ acts continuously by
isometries. This is widely used in geometric analysis. It is also custom to consider the free loop space
of continuous closed curves $LM=C^0(S^1,M)$ in topology. Since $\Lm M$ and $LM$ are $G$-weakly homotopy
equivalent with any subgroup $G$ of $O(2)$ (cf. Hingston \cite{Hingston1984}, p.101), we shall use both
notations and do not distinguish them in this paper.

It is well known (cf. Chapter 1 of \cite{Kli1978}) that $\gamma$ is a closed geodesic or a constant curve
on $(M,F)$ if and only if $\gamma$ is a critical point of the energy functional
$$ E(c)=\frac{1}{2}\int_{0}^{1}F(\dot{c},\dot{c})dt. $$

For more studies on the closed geodesics, we refer the readers to the survey papers of Bangert
\cite{Bangert1985}, Long \cite{Lo2006}, and Taimanov \cite{Tai2010}. Among others, two important topological
invariants associated to the free loop space $\Lm M$ of a compact manifold $M$ are used by many
mathematicians in the study of multiplicity and stability of closed geodesics on $M$.

The first one is the Betti number sequence $\{\beta_k(\Lm M;\mathbb{F})\}_{k\in\Z}$ of the free loop space
$\Lm M$ with $\beta_k(\Lm M;\mathbb{F}) = \rank H_{*}(\Lm M;\mathbb{F})$ for an arbitrary field $\mathbb{F}$.
In 1969 Gromoll and Meyer (\cite{GM1969JDG}, Theorem 4) established  the existence of infinitely many distinct
closed geodesics on $M$, provided that $\{{\beta}_k(\Lm M;\mathbb{Q})\}_{k\in\Z}$ is unbounded. Then
Vigu$\acute{e}$-Poirrier and Sullivan \cite{VS1976} further proved in 1976 that for a compact simply connected manifold
$M$, the Gromoll-Meyer condition holds if and only if $H^{*}(M;\mathbb{Q})$ is generated by more than one element.
Although the Gromoll-Meyer theorem is valid actually for any field $\mathbb{F}$, and there are really some spaces
with bounded $\{{\beta}_k(\Lm M;\mathbb{Q})\}_{k\in\Z}$ but unbounded $\{{\beta}_k(\Lm M;\mathbb{Z}_{2})\}_{k\in\Z}$,
it can not be applied to the compact globally symmetric spaces of rank $1$ which consist in
\be   S^{n},\ \mathbb{R}P^{n},\ \mathbb{C}P^{n},\ \mathbb{H}P^{n}\ \text{and}\ {\rm CaP}^{2},  \label{mflds}\ee
since $\{{\beta}_k(\Lm M,\mathbb{F})\}_{k\in\Z}$ with $M$ in (\ref{mflds}) is bounded with respect to any field
$\mathbb{F}$ (cf. Remark in \cite{Ziller1977}, pp. 17-18 of Ziller).

The second one is the Betti number sequence $\{\bar{\beta}_k(\Lm M,\Lm^{0} M;\mathbb{F})\}_{k\in\Z}$ of the
free loop space pair $(\Lm M,\Lm^{0} M)$ with $\bar{\beta}_k(\Lm M,\Lm^{0} M;\mathbb{F})$ being the rank of
the relative $S^1$-equivariant homology $H_{*}^{S^1}(\Lm M,\Lm^{0}M; \mathbb{F})$ defined by
\be  H_{*}^{S^1}(\Lm M,\Lm^{0}M; \mathbb{F})
      = H_{*}(\Lm M\times_{S^1} ES^{1},\Lm^{0}M\times_{S^1} ES^{1}; \mathbb{F}), \lb{S1homology}\ee
where $\Lm^{0}M=M$ is the set of constant curves on $M$. When $\mathbb{F}=\mathbb{Q}$, it can be proved further that
\be  H_{*}^{S^1}(\Lm M,\Lm^{0}M; \mathbb{Q}) \cong H_{*}(\Lm M/{S^1},\Lm^{0}M/{S^1}; \mathbb{Q}),  \label{bridge1}\ee
which enables Rademacher (\cite{Rad1989} of 1989 and \cite{Rad1992} of 1992) to establish the resonance identity
of prime closed geodesics on the simply connected manifolds of (\ref{mflds}) provided that their total number is
finite. This identity relates local topological invariants of prime closed geodesics to the average Betti number
$\bar{B}(\Lm M, \Lm^{0}M; \mathbb{Q})$ defined by
\be  \bar{B}(\Lm M, \Lm^{0}M; \mathbb{Q})
  = \lim_{m\to+\infty}\frac{1}{m}\sum_{k=0}^{m}(-1)^k\bar{\beta}_k(\Lm M,\Lm^{0}M; \mathbb{Q}), \lb{aBLM}\ee
and then was used by many authors to study the multiplicity of closed geodesics on such manifolds. One example is
the existence of at least two distinct closed geodesics on every $2$-dimensional Finsler sphere $S^{2}$ proved by
Bangert and Long (\cite{BL2010} of 2010).

But for the multiplicity of closed geodesics on non-simply connected manifolds whose free loop space possesses
bounded Betti number sequence, we are aware of not many works. For example, in 1981, Ballman, Thorbergsson
and Ziller \cite{BTZ1981} proved that every Riemannian manifold with the fundamental group being a nontrivial
finitely cyclic group and possessing a generic metric has infinitely many distinct closed geodesics. In 1984,
Bangert and Hingston \cite{BH1984} proved that any Riemannian manifold with fundamental group an infinite
cyclic group has infinitely many distinct closed geodesics. To use Morse theory to study this problem, one should
know the global topological structure of the free loop space on these manifolds. As far as the authors know, there
seems to be only two precise results on real projective spaces obtained by Westerland \cite{West2005}, \cite{West2007}
in the field $\Z_2$.

However when one tries to apply Westerland's results directly to study the multiplicity of closed geodesics on
$\mathbb{R}P^{d}$, two difficulties appear which also lie in the simply connected cases if the coefficient
field used is $\Z_{2}$. The first one is that for every $M$ in (\ref{mflds}), the following isomorphism
\be   H_{*}^{S^1}(\Lm M,\Lm^{0}M; \mathbb{Z}_{2}) \cong H_{*}(\Lm M/{S^1},\Lm^{0}M/{S^1}; \mathbb{Z}_{2}),
       \label{bridge2}\ee
may not hold, which would however be crucial to relate the relative $S^1$-equivariant homogloy
$H_{*}^{S^1}(\Lm M,\Lm^{0}M;\Z_{2})$ to the $\Z_2$-critical modules of closed geodesics (cf. the definition in
(\ref{CritG})). Note that such an isomorphism (\ref{bridge3}) is proved in our Lemma \ref{Lm3.3} below for the
non-contractible free loop space.
The second one is that the $\Z_2$-critical module of even $m$-iterates of a closed geodesic $c$ is related to
the parity of the difference $i(c^m)-i(c)$ of Morse indices and may lead to unbounded Morse type number sequence
so that Gromoll-Meyer type argument is not applicable. We refer readers to Lemma 4.1.4 (ii) on p.127 of
\cite{Kli1978} as well as Proposition 3.8 and its proof in pp. 345-346 of \cite{BL2010} for more details.

In this paper we overcome the above two difficulties for $\mathbb{R}P^{2n+1}$ by restricting the problem of closed
geodesics to the non-contractible component of the free loop space of $\RP$. More precisely, let $M=\mathbb{R}P^{d}$ with $d\ge 2$.
Then $\pi_{1}(M) = \mathbb{Z}_{2} = \{e, g\}$ with $e$ being the identity and $g$ being the generator of
$\mathbb{Z}_{2}$ satisfying $g^{2}=e$ and the free loop space $LM$ possesses a natural decomposition
$$    LM = L_{e}M\bigsqcup L_{g}M,   $$
where $L_{e}M$ and $L_{g}M$ are the two connected components of $LM$ whose elements are homotopic to $e$ and $g$
respectively. We shall prove in Lemma \ref{Lm3.3} below
\begin{equation}
\label{bridge3}
H_{*}^{S^1}(L_{g}M;\mathbb{Z}_{2})\cong H_{*}(L_{g}M/S^{1};\mathbb{Z}_{2}).
\end{equation}
Moreover, we observe that every closed geodesic on $\mathbb{R}P^{d}$ is orientable if and only if
$d\in 2\N-1$, and an $m$-th iterate $c^m$ of a non-contractible closed geodesic $c$ on $\mathbb{R}P^{2n+1}$ is
still non-contractible if and only if $m$ is odd, which then implies $i(c^m)-i(c)$ for odd $m$ is always even.
This makes the structure of the $\Z_2$-critical modules of odd iterates of non-contractible closed geodesics
be similar to the case in the coefficient field $\mathbb{Q}$ and become rather simple. Thus it is possible to
get some information on them, if we can get the topological structure of $L_{g}\RP$, which is our first goal in
this paper. To this end, we use the following ideas.

(i) Consider the fibrations
$$  L_{e}M\to L_{e}M\times_{S^{1}}ES^{1}\to BS^{1},  $$
and
$$¡¡L_{g}M\to L_{g}M\times_{S^{1}}ES^{1}\to BS^{1},  $$
as well as their Leray-Serre spectral sequences which we will denote for simplicity by $\LS(L_eM)$ and
$\LS(L_gM)$ respectively in this paper.

Define the $S^1$-equivariant Betti numbers and the Poincar\'e series of $L_{*}M$ with $L_{*}M = LM$, $L_{e}M$
or $L_{g}M$  by
\bea
\bar{\bb}_k(L_{*}M; \Z_{2}) &=& \dim H_k^{S^1}(L_{*}M;\Z_2),   \lb{BettiN}\\
P^{S^1}(L_{*}M;\Z_2)(t) & =& \sum_{k=0}^{\infty}\bar{\bb}_k(L_{\ast}M; \Z_{2})t^k. \lb{Po1}
\eea
Here we do not use the relative forms due to the fact that $\Lm^{0} M\cap L_{g}M=\emptyset.$ It then follows
\be P^{S^1}(LM;\mathbb{Z}_{2})(t) = P^{S^1}(L_{e}M;\mathbb{Z}_{2})(t)+P^{S^1}(L_{g}M;\mathbb{Z}_{2})(t). \lb{Po2}\ee

In Section 4 of this paper, we apply Chas-Sullivan theory and Leray-Serre spectral sequence to carry out related
computations to obtain in Theorem \ref{ThmBVstr} below four possible Batalin-Vilkovisky (we write B-V below for
short) algebraic structures of $\mathbb{H}_{*}(LM;\mathbb{Z}_2)=H_{*+(2n+1)}(LM;\mathbb{Z}_2)$ according to the behaviors
of its generators $\td{x}$, $\td{v}$ and $\td{w}$ constructed in Lemma \ref{algebraic1}.

(ii) In \cite{Katok1973} of 1973, Katok constructed a famous family of Finsler metrics $N_{\aa}$ for
$\aa\in (0,1)\bs\Q$ on spheres, specially $S^{2n+1}$, which possesses precisely $2n+2$ distinct prime closed
geodesics. In Section 3 below, via inducing the metrics $N_{\aa}$ to $\mathbb{R}P^{2n+1}$, we then prove that
the $S^1$-equivariant Betti number sequence $\{\bar{\bb}_k(L_gM; \Z_{2})\}_{k\in\Z}$ is bounded via the
boundedness of the Morse number sequence. Such a result will help us simplify the proof of Theorem \ref{Thm1.1}.

(iii) In Section 5 of this paper, we compute the $S^1$-equivariant Poincar\'e series associated to the third pages
of $\LS(L_eM)$ and $\LS(L_gM)$ under each of the four possible B-V structures, and then comparing  their sums
with $P^{S^1}(LM;\Z_{2})(t)$ computed by Westerland in \cite{West2007}, we obtain the $S^1$-equivariant Poincar\'e
series of $L_gM$, which then yields  the average $S^1$-equivariant Betti number of $L_gM$  defined by
\be \bar{B}(L_gM ;\Z_{2}) = \lim_{m\to+\infty}\frac{1}{m}\sum_{k=0}^{m}(-1)^k\bar{\beta}_k(L_gM ;\Z_{2}). \lb{aB.1}\ee
The following is our first main result in this paper.

\begin{theorem}\label{Thm1.1}
For $M=\mathbb{R}P^{2n+1}$ with $n\ge 1$, $\LS(L_eM)$ and $\LS(L_gM)$ collapse at the
second page and the third one respectively. Moreover, the $S^1$-equivariant Poincar\'e series of $L_gM$ satisfies
\be P^{S^1}(L_{g}M;\mathbb{Z}_{2})(t) =  \frac{1-t^{2n+2}}{(1-t^{2n})(1-t^{2})}, \lb{Clm1.1}\ee
and the average $S^1$-equivariant Betti number of $L_gM$ satisfies
\be  \bar{B}(L_gM;\Z_{2}) = \frac{n+1}{2n}.   \lb{aB.2}\ee
\end{theorem}

\begin{remark}
For the complex projective space $\mathbb{C}P^{n}$, B\"{o}kstedt and Ottosen \cite{BO2007} proved in 2007
that $\LS(L\mathbb{C}P^{n})$ collapses at the third page.
\end{remark}

Next we apply Theorem \ref{Thm1.1} to study the problem of non-contractible closed geodesics on the Finsler
manifold $M=(\mathbb{R}P^{2n+1},F)$.

Recall that on a Finsler manifold $(M,F)$, for some integer $m\in\N$, the $m$-th iterate $c^m$ of $c\in \Lm M$ is
defined by
$$  c^m(t) = c(mt), \qquad \forall\;t\in [0,1].   $$
The inverse curve $c^{-1}$ of $c$ is defined by $c^{-1}(t)=c(1-t)$ for $t\in S^1$.

\begin{definition}\label{def-CG(A)}
Let $(M,F)$ be a Finsler manifold (or a Riemannian manifold), and $A$ be a subset of the free loop space $\Lm M$.
A closed curve $c:S^1=\R/\Z\to M$ belonging to $A$ is {\it prime} in $A$, if it is not a multiple covering (i.e.,
iterate) of any other closed curve on $M$ belonging to $A$. Two prime closed curves $c_1$ and $c_2$ in $A$ are
{\it distinct} (or {\it geometrically distinct}), if they do not differ by an $S^1$-action (or $O(2)$-action).
We denote by $\CG(A)$ the set of all closed geodesics on $(M,F)$ (or in Riemannian case) which are distinct
prime in $A$. As usual we denote by $\CG(M,F)=\CG(\Lm M)$ when $A=\Lm M$.
\end{definition}

For the Finsler manifold $M=(\mathbb{R}P^{2n+1},F)$, as in \cite{BL2010}, let $c$ be a closed geodesic in $M$
satisfying the following isolated condition:
$$ S^{1}\cdot c^{m}\ \text{is\ an\ isolated\ critical\ orbit\ of}\ E \ \text{for\ every}\ m\ge 1. \ \leqno{({\rm ISO})} $$
We let
$$  \Lambda_{g}(c) = \{\gamma\in \Lambda_{g}M\mid E(\gamma)<E(c)\},  $$
and define the $\mathbb{Z}_{2}$-critical module of $c^{2m-1}$ by
\be \bar{C_{q}}(E,c^{2m-1};\Z_{2})
   = H_{q}\left((\Lambda_{g}(c^{2m-1})\cup S^{1}\cdot c^{2m-1})/S^{1},\Lambda_{g}(c^{2m-1})/S^{1};\Z_{2}\right).  \lb{CritG}\ee
As we shall prove in Section 3 below, we have
\be \bar{C_{q}}(E,c^{2m-1})
= H_{q-i(c^{2m-1})}(N_{c^{2m-1}}^{-}\cup\{c^{2m-1}\},N_{c^{2m-1}}^{-};\mathbb{Z}_{2}), \lb{shift.0}\ee
where $N_{c^{2m-1}}^{-}=N_{c^{2m-1}}\cap\Lambda_{g}(c^{2m-1})$, $N_{c^{2m-1}}$ is the local characteristic
manifold of $E$ at $c^{2m-1}$, and to get (\ref{shift.0}) we have used the fact $i(c^{2m-1})-i(c)\in 2\Z$
proved in (\ref{Bott.1}) below. Here properties of odd iterates of $c$ are crucial.

As usual, for $m\in\N$ and $l\in\Z$ we define the local homological type numbers of $c^{2m-1}$ by
\be k_{l}(c^{2m-1})
= \dim H_{l}(N_{c^{2m-1}}^{-}\cup\{c^{2m-1}\},N_{c^{2m-1}}^{-};\mathbb{Z}_{2}).  \lb{CGht1}\ee

Based on works of Rademacher in \cite{Rad1989}, Long and Duan in \cite{LD2009} and \cite{DL2010},
we define the {\it analytical period} $n_c$ of the closed geodesic $c$ by
\be n_c = \min\{j\in 2\N\,|\,\nu(c^j)=\max_{m\ge 1}\nu(c^m),\;\;
                  \forall\,m\in 2\N-1\}. \lb{CGap1}\ee
Note that here in order to simplify the study for non-contractible closed geodesics in $\mathbb{R}P^{2n+1}$,
we have slightly modified the definition in \cite{LD2009} and \cite{DL2010} by requiring the analytical
period to be even. Then by the same proofs in \cite{LD2009} and \cite{DL2010}, we have
\be  k_{l}(c^{2m-1+hn_c}) = k_{l}(c^{2m-1}), \qquad \forall\;m,\;h\in \N,\;l\in\Z.  \lb{CGap2}\ee
For more detailed properties of the analytical period $n_c$ of a closed geodesic $c$, we refer readers to
the two Section 3s in \cite{LD2009} and \cite{DL2010}.

As in \cite{BK1983}, we have

\begin{definition}\label{def-hv} Let $(M,F)$ be a compact Finsler manifold. A closed geodesic $c$ on $M$
is homologically visible, if there exists an integer $k\in\Z$ such that $\bar{C}_k(E,c) \not= 0$. We denote
by $\CG_{\hv}(M,F)$ the set of all distinct homologically visible prime closed geodesics on $(M,F)$.
\end{definition}

It is well known that on a compact Finsler manifold $M$, there exists at least one homologically visible prime
closed geodesic, because the topology of the free loop space $\Lm M$ is non-trivial (cf. \cite{Bangert1985}).
Motivated by the resonance identity proved in \cite{Rad1989}, in Section 6 below, as an application of Theorem
\ref{Thm1.1} we obtain the following resonance identity on the non-contractible closed geodesics on Finsler
$M=(\mathbb{R}P^{2n+1},F)$. Note that here if $\;^{\#}\CG(M)<+\infty$ and $c$ is a prime homologically
visible closed geodesic on $M$, then $\hat{i}(c) > 0$ must hold by Lemma \ref{Lm3.4} below.

\begin{theorem}\label{Thm1.2} Suppose the Finsler manifold $M=(\mathbb{R}P^{2n+1},F)$ possesses only finitely
many distinct prime closed geodesics, among which we denote the distinct non-contractible homologically visible prime closed
geodesics by $c_1, \ldots, c_r$  for some integer $r>0$. Then we have
\be  \sum_{j=1}^{r}\frac{\hat{\chi}(c_j)}{\hat{i}(c_j)} = \bar{B}(\Lm_gM ;\Z_{2}) = \frac{n+1}{2n},   \lb{reident1}\ee
where the mean Euler number $\hat{\chi}(c_j)$ of $c_j$ is defined by
$$  \hat{\chi}(c_j) = \frac{1}{n_j}\sum_{m=1}^{n_j/2}\sum_{l=0}^{4n}(-1)^{l+i(c_{j})}k_{l}(c_{j}^{2m-1}), $$
and $n_j=n_{c_j}$ is the analytical period of $c_j$.
\end{theorem}

\begin{remark} (i) Note that homologically invisible closed geodesics, if they exist, have no contributions
to the resonance identity (\ref{reident1}).

(ii) For the special case when each $c_{j}^{2m-1}$ is non-degenerate with $1\leq j\leq r$ and $m\in\mathbb{N}$,
we have $n_{j}=2$ and $k_l(c_j)=1$ when $l=0$, and $k_l(c_j)=0$ for all other
$l\in\Z$. Then (\ref{reident1}) has the following simple form
\be  \sum_{j=1}^{r}(-1)^{i(c_{j})}\frac{1}{\hat{i}(c_{j})}=\frac{n+1}{n}.  \lb{reident2}\ee
\end{remark}

\setcounter{figure}{0}
\setcounter{equation}{0}
\section{Preliminaries}%Section 2

In this section, we recall some facts on the free loop spaces and Leray-Serre spectral sequence.

\subsection{String topology on the loop homology}

In their seminar paper \cite{CS1999}, Chas and Sullivan introduced a new collection of invariants of manifolds for
the loop homology. Specifically, let $M$ be a $d$-dimensional oriented compact manifold and $LM$ be the free loop
space of $M$. They geometrically constructed the loop product $\bullet$ and the loop bracket $\{,\}$ on the loop
homology $\mathbb{H}_{*}(LM)=H_{*+d}(LM),$  and proved that $(\mathbb{H}_{*}(LM),\bullet,\{,\})$ is a Gerstenhaber
algebra.
\begin{definition} $(A,\bullet,\{,\})$ is called a Gerstenhaber algebra, if

$(1)$ $(A,\bullet)$ is a graded commutative, associative algebra,

$(2)$ $\{,\}$ is a Lie bracket of degree $+1$, that is , for every $a,b,c\in A$,

\qquad $(i)$ $\{a,b\}=-(-1)^{(|a|+1)(|b|+1)}\{b,a\},$

\qquad $(ii)$ $\{a,\{b,c\}\}=\{\{a,b\},c\}+(-1)^{(|a|+1)(|b|+1)}\{b,\{a,c\}\}$,

$(3)$ $\{a, b\bullet c\}=\{a,b\}\bullet c+(-1)^{(|a|+1)|b|}b\bullet\{a,c\},$\\
where $|a|$ denotes the degree of $a$.
\end{definition}

Consider the circle action $\eta:S^{1}\times LM\to LM$, defined by
$$\eta(\theta,\gamma)(t)=\gamma(\theta+t),\ \forall (\theta,\gamma)\in S^{1}\times LM.$$ Then,
$\eta$ induces a degree $+1$ operator $\bigtriangleup:H_{*}(LM)\to H_{*+1}(LM)$ defined by
\begin{equation}
\label{CSoperator}
\bigtriangleup(v)=\eta_{*}([S^{1}]\otimes v),\qquad \forall v\in H_{*}(LM),
\end{equation}
with $[S^{1}]$ the generator of $H_{1}(S^{1})$. In \cite{CS1999}, the authors also proved that
$(\mathbb{H}_{*}(LM),\bullet,\bigtriangleup)$ is a Batalin-Vilkovisky algebra.
\begin{definition} $(A,\bullet,\bigtriangleup)$ is called a Batalin-Vilkovisky algebra, if

$(1)$ $(A,\bullet)$ is a graded commutative, associative algebra,

$(2)$  $\bigtriangleup\circ\bigtriangleup=0,$

$(3)$ $(-1)^{|a|}\bigtriangleup(a\bullet b)-(-1)^{|a|}\bigtriangleup(a)\bullet b-a\bullet \bigtriangleup(b)$ is
a derivation of each variable.
\end{definition}
In particular for the free loop space $LM$, we still have the following B-V formulae
\begin{equation}
\label{BV1}
\{a,b\}=(-1)^{|a|}\bigtriangleup(a\bullet b)-(-1)^{|a|}\bigtriangleup(a)\bullet b-a\bullet \bigtriangleup(b),
\end{equation} (see Corollary 5.3 in \cite{CS1999}). Since the coefficient filed used in this paper is $\mathbb{Z}_{2}$,
the B-V formulae (\ref{BV1}) has the simpler form
\begin{equation}
\label{leibnizlike}\bigtriangleup(a\bullet b)=\bigtriangleup(a)\bullet b+a\bullet \bigtriangleup(b)+\{a,b\},
\end{equation}
which is similar as but a little more complicated than the Leibniz formulae.

In 2002, Cohen and Jones \cite{CJ2002} realized Chas-Sullivan loop product $\bullet$ as the cup product $\cup$ in
the Hochschild cohomology $HH^{*}(H^{*}(M), H^{*}(M))$, when $M$ is simply connected, that is, they succeeded in
establishing a ring isomorphism
$$(\mathbb{H_*}(LM),\bullet)\stackrel{\cong}{\longrightarrow}(HH^{*}(H^{*}(M), H^{*}(M)),\cup),$$
(see Theorem 3 in \cite{CJ2002} for details). Such a result also holds for $\mathbb{R}P^{d}$ due to Lemma 5.4 in
\cite{West2005}.

\subsection{Leray-Serre spectral sequence}

The theory of Leray-Serre spectral sequence can be found in many literatures such as Hatcher \cite{Ha} and McCleary
\cite{Mc2001}. We sketch it for the reader's convenience.

Consider the fibration
$$ F\to E\to B, $$
where $F$ is the connected fiber, $E$ is the total space and the base space $B$ is simply connected. Since the
coefficient field in this paper is $\mathbb{Z}_{2}$, there is a cohomology Leray-Serre spectral sequence
$\{E_{r}^{p,q}, \hat{d}_{r}\}$ converging to $H^{*}(E)$. Moreover, the second page
$$E_{2}^{p,q}\thickapprox H^{p}(B; H^{q}(F))\thickapprox H^{p}(B)\otimes H^{q}(F).$$
Similarly, we have the homology Leray-Serre spectra sequence $\{E_{p,q}^{r}, {d}_{r}\}$ converging to $H_{*}(E)$.
Generally speaking, the cohomology form is more powerful than the homology one, since the cohomology groups equipped
with the cup product becomes a cohomology ring. By the Leibniz formulae, the differential $\hat{d}_{r}$ is determined
if one has known the value of $\hat{d}_{r}$ on the generators of the ring.

However, we use in this paper  the homology Leray-Serre spectral sequence based on the following observations. Due to
the work of \cite{CJ2002} and \cite{West2005}, the ring structure
$(\mathbb{H}_{*}(L\mathbb{R}P^{2n+1};\mathbb{Z}_{2}),\bullet)$ is clear. Since
$(\mathbb{H}_{*}(L\mathbb{R}P^{2n+1};\mathbb{Z}_{2}),\bullet,\bigtriangleup)$ is a B-V algebra, the close relationship
between the operator $\bigtriangleup$  with the differential $d_{2}$  enables us to compute the third page $E_{p,q}^{3}$
by use of the B-V formulae instead of the Leibniz formulae in the cohomology case.

Let $Y$ be a connected $S^{1}$-space  with the action map $\eta:S^{1}\times Y\to Y.$ Consider the composition map
$$ Y\stackrel{\tau}{\longrightarrow}S^{1}\times Y\stackrel{\eta}{\longrightarrow}Y, $$
where $\tau(y)=(1,y).$ Since $\eta\circ\tau=id$, the map
$$\eta^{*}:H^{*}(Y;\mathbb{Z}_{2})\to H^{*}(S^{1}\times Y;\mathbb{Z}_{2}),$$
is injective. Denote by $\{S^1\}$ the generator of $H^{1}(S^{1};\mathbb{Z}_{2})$, then $\eta^{*}$ induces a degree $-1$
map $d:H^{*}(Y;\mathbb{Z}_{2})\to H^{*-1}(Y;\mathbb{Z}_{2})$ satisfying that
\begin{equation}
\label{924b}
\eta^{*}(y)=1\otimes y+\{S^{1}\}\otimes dy.
\end{equation}
It is known that $d$ is a linear differential operator, i.e., $d\circ d=0$, $d(x+y)=dx+dy$
and $d(xy)=xdy+ydx$ (see \cite{BO1999} Proposition 3.2, p. 253).

\begin{lemma} (\cite{BO1999} Proposition 3.3.)
\label{bo1} The fibration $$Y\to Y\times_{S^1}ES^{1}\to BS^{1}$$
has the following cohomology  Leray-Serre spectral sequence:
$$E_{2}^{*,*}=H^{*}(BS^{1})\otimes H^{*}(Y)\Rightarrow
H^{*}(Y\times_{S^1}ES^{1}),$$ The differential in the $E_2$-page is given by $\hat{d}_{2}:H^{*}(Y)\to \hat{u}H^{*}(Y)$,
$$\hat{d}_{2}(1\otimes y)=\hat{u}\otimes d(y)£¬$$ where $d$ is the  differential defined in (\ref{924b}) and $\hat{u}$
is the generator of $H^{*}(BS^1)$ with degree $2$.
\end{lemma}

\begin{remark}
\label{denglingyu0} The above lemma has a more general form:
$$\hat{d}_{2}(\hat{u}^{p}\otimes y)=\hat{u}^{p+1}\otimes d(y),\
\forall\ p\in\mathbb{N}\cup\{0\}\ \text{and}\ y\in H^{*}(Y).$$
Actually, since $\hat{u}^{p}\otimes1$ lies in the horizontal axis of the first quadrant,
we have $\hat{d}_{2}(\hat{u}^{p}\otimes 1)=0$ for $p\in\mathbb{N}\cup\{0\}$. By Leibniz formulae, we obtain
$$\begin{aligned}
\hat{d}_{2}(\hat{u}^{p}\otimes y)&=\hat{d}_{2}((\hat{u}^{p}\otimes1)\cup(1\otimes
y))\\
&=\hat{d}_{2}(\hat{u}^{p}\otimes1)\cup(1\otimes
y)+(\hat{u}^{p}\otimes1)\cup \hat{d}_{2}(1\otimes y)\\
&=(\hat{u}^{p}\otimes1)\cup \hat{d}_{2}(1\otimes y)\\
&=(\hat{u}^{p}\otimes1)\cup  (\hat{u}\otimes y)\\
&=\hat{u}^{p+1}\otimes y. \end{aligned}$$
\end{remark}
In the sequel, we denote $\hat{a}\in H^{*}(Y)$ the dual of $a$ for every $a\in H_{*}(Y)$, then $\hat{a}(a)=\langle\hat{a},a\rangle=1.$
Corresponding to Lemma \ref{bo1}, we have

\begin{lemma}\label{xym1}
The fibration $Y\to Y\times_{S^1}ES^{1}\to BS^{1}$ has the following homology Leray-Serre spectral sequence:
$$E_{*,*}^{2}=H_{*}(BS^{1})\otimes H_{*}(Y)\Rightarrow H_{*}(Y\times_{S^1}ES^{1}),$$ the differential $d_2$ in the $E^2$-page is given by
${d}_{2}:u^{p}\otimes H_{*}(Y)\to u^{p-1}\otimes H_{*}(Y)$ Ϊ
$${d}_{2}(u^{p}\otimes v)=u^{p-1}\otimes\bigtriangleup(v),\ \forall\ p\in\mathbb{N},$$ where $u\in H_{2}(BS^1)$
is the dual of $\hat{u}$ and $\bigtriangleup(v)=\eta_{*}([S^{1}]\otimes v)$ with $[S^{1}]$ the dual of $\{S^{1}\}$.
\end{lemma}
\Proof
On the one hand, for every $v\in H_{q-1}(Y)$ and $y\in H^{q}(Y)$ we have  from $\eta:S^{1}\times Y {\longrightarrow}Y$ that
\begin{equation*}
\begin{aligned}
\langle\eta^{*}y,[S^{1}]\otimes v\rangle&=\langle1\otimes y+\{S^{1}\}\otimes dy,[S^{1}]\otimes v\rangle\\
&=\langle\{S^{1}\}\otimes dy,[S^{1}]\otimes v\rangle\\
&=\langle dy, v\rangle.
\end{aligned}
\end{equation*}
On the other hand,
\begin{equation*}
\begin{aligned}
\langle\eta^{*}y,[S^{1}]\otimes v\rangle= \langle
y,\eta_{*}([S^{1}]\otimes v)\rangle=\langle
y,\bigtriangleup(v)\rangle.
\end{aligned}
\end{equation*}
So,
\begin{equation}
 \label{2012fangbian}
 \langle dy,v\rangle=\langle y,\bigtriangleup(v)\rangle.
\end{equation}¡¡

Let $w=u^{p}\otimes v$, then by the duality of cohomology and homology Leray-Serre spectral sequences,
Remark \ref{denglingyu0} and (\ref{2012fangbian}),¡¡
\begin{equation*}
\begin{aligned}
\langle \hat{u}^{p-1}\otimes y, d_{2}w\rangle&=\langle
\hat{d}_{2}(\hat{u}^{p-1}\otimes y),w\rangle\\
&=\langle
\hat{u}^{p}\otimes dy,w\rangle\\
&=\langle\hat{u}^{p}\otimes dy,u^{p}\otimes v\rangle\\
&=\langle dy,v\rangle\\
&=\langle y, \bigtriangleup(v)\rangle\\
&=\langle
\hat{u}^{p-1}\otimes y,
u^{p-1}\otimes\bigtriangleup(v)\rangle,
\end{aligned}
\end{equation*}
 that is, $d_{2}(u^{p}\otimes
v)=u^{p-1}\otimes\bigtriangleup(v)$. $\hfill\Box$

\begin{remark}
In the sequel, we actually do computations in $\{\mathbb{E}^{2},d_{2}\}$, where
$\mathbb{E}^{2}=H_{*}(BS^{1})\otimes \mathbb{H}_{*}({Y}),$ that is, we only shift the degree of $H_{*}(Y)$ while
do not shift those of $H_{*}(BS^{1})$. Since it is only the change of notation, the above lemma is still valid.
\end{remark}

\setcounter{figure}{0}
\setcounter{equation}{0}
\section{Katok's metrics on spheres and real projective spaces}%Section 3

In this section, we prove that the Betti number sequence of $P^{S^{1}}(L_{g}\mathbb{R}P^{2n+1};\mathbb{Z}_{2})$
is  bounded via Katok's famous metrics on $S^{2n+1}$. Such
a result will help us to simplify the proof of Theorem \ref{Thm1.1}.

In 1973, Katok \cite{Katok1973} constructed his famous irreversible Finsler metrics on $S^{n}$ which possess
only finitely many distinct prime closed geodesics. His examples were further studied closely by Ziller
\cite{Ziller1982} in 1982, from which we borrow most of the notations for the particular case $S^{2n+1}$.

Let $S^{2n+1}$ be the standard sphere with the canonical Riemannian metric $g$ and the one-parameter group
of isometries
$$   \phi_{t}=\diag(R(pt/p_{1}),\dots,R(pt/p_{n+1})),   $$
where $p_{i}\in\mathbb{Z}$, $p=p_{1}\cdots p_{n+1}$ and $R(\omega)$ is a rotation in $\mathbb{R}^{2}$ with
angle $\omega$. Let $TS^{2n+1}$ and $T^{*}S^{2n+1}$ be its tangent bundle and cotangent bundle respectively.
Define $H_{0}, H_{1}:T^{*}S^{2n+1}\to\mathbb{R}$ by
$$ H_{0}(x)=\|x\|_{*}\quad \text{and}\quad H_{1}(x)=x(V),\quad \forall x\in T^{*}S^{2n+1},  $$
where $\|\cdot\|_{*}$ denotes the dual norm of $g$ and $V$ is the vector field generated by $\phi_{t}$. Let
$$  H_{\alpha}=H_{0}+\alpha H_{1} \qquad {\rm for}\quad \alpha \in (0,1).  $$
Then $\frac{1}{2}H_{\alpha}^{2}$ is homogeneous of degree two and the Legendre transform
$$  L_{\frac{1}{2}H_{\alpha}^{2}}=D_{F}\left(\frac{1}{2}H_{\alpha}^{2}\right):T^{*}S^{2n+1}\to TS^{2n+1},  $$
is a global diffeomorphism. Hence,
$$    N_{\alpha}=H_{\alpha}\circ L_{\frac{1}{2}H_{\alpha}^{2}}^{-1}   $$
defines a Finsler metric on $S^{2n+1}$. Since $H_{\alpha}(-x)\neq H_{\alpha}(x)$, $N_{\alpha}$ is not reversible.
It was proved that $(S^{2n+1}, N_{\alpha})$ with $\alpha \in (0,1)\bs\Q$ possesses precisely $2(n+1)$ distinct
prime closed geodesics (cf. Katok \cite{Katok1973} and pp. 137-139 of Ziller \cite{Ziller1982} for more details).

Consider the antipodal map $A:S^{2n+1}\to S^{2n+1}$ defined by $ A(p)=-p.$ We first prove

\begin{lemma}\label{Lm3.1}
$A:(S^{2n+1}, N_{\alpha})\to (S^{2n+1}, N_{\alpha})$ is an isometry.
\end{lemma}

\Proof  For any $p\in S^{2n+1}$, we consider the following diagram
\be
\begin{tabular}{ccc}
$T_p S^{2n+1}$                           & $\mapright{A_{\ast}}$ & $T_{-p}S^{2n+1}$ \\
$\mapdown{L_{\frac{1}{2}H_{\aa}^2}^{-1}}$&                       & $\mapdown{L_{\frac{1}{2}H_{\aa}^2}^{-1}}$ \\
$T_p^{\ast}S^{2n+1}$                     & $\mapright{A^{\ast}}$ & $T_{-p}^{\ast}S^{2n+1}$ \\
$\qquad\mapse{H_{\aa}}$                  &                       & $\mapsw{H_{\aa}}\qquad$ \\
                                         & $\mathbb{R}$                  &  \\
\end{tabular}  \lb{1}\ee

{\bf Claim 1.} {\it The lower triangle in (\ref{1}) commutes.}

In fact, for any $x\in T_{p}^{*}S^{2n+1},$  we have $A^{*}(x)=-x\in T_{-p}^{*}S^{2n+1}$ and
\bea  H_{\alpha}\circ A^{*}(x)
&=& H_{\alpha}(-x)   \nn\\
&=& \|-x\|_{*}+(-x)(V_{-p})  \nn\\
&=& \|x\|_{*}+(-x)(-V_{p})  \nn\\
&=& \|x\|_{*}+x(V_{p})  \nn\\
&=& H_{\alpha}(x),   \label{1217c}\eea
where the third identity is due to the fact that $A$ is isomorphic for the canonical metric $\|\cdot\|_{*}$ and
$V_{-p}=-V_{p}$ which follows by the definition of $V$. This proves the claim.

To prove that $A$ is an isometry of $S^{2n+1}$, i.e.,
\begin{equation}\label{1217a}
H_{\alpha}\circ L_{\frac{1}{2}H_{\alpha}^{2}}^{-1}\circ A_{*}(X)
           = H_{\alpha}\circ L_{\frac{1}{2}H_{\alpha}^{2}}^{-1}(X),\ \forall X\in T_{p} S^{2n+1},
\end{equation}
it is sufficient by (\ref{1217c}) to prove

{\bf Claim 2.} {\it The upper square in (\ref{1}) commutes, i.e.,
\begin{equation}\label{1217b}
L_{\frac{1}{2}H_{\alpha}^{2}}^{-1}\circ A_{*}(X)
  =  A^{*}\circ L_{\frac{1}{2}H_{\alpha}^{2}}^{-1}(X),\ \forall X\in T_{p} S^{2n+1},
\end{equation}}

In fact, (\ref{1217b}) is equivalent to
\begin{equation}\label{1217d}
A_{*}\circ L_{\frac{1}{2}H_{\alpha}^{2}}(x)=L_{\frac{1}{2}H_{\alpha}^{2}}\circ A^{*}(x),\ \forall x\in T_{p}^{*} S^{2n+1}.
\end{equation}

Note that for any $x\in T_{p}^{*}S^{2n+1} $,
$$\begin{aligned} L_{\frac{1}{2}H_{\alpha}^{2}}(x)
&=D_{F}\left(\frac{1}{2}H_{\alpha}^{2}\right)(x)\\
&=H_{\alpha}(x)\cdot D_{F}H_{\alpha}(x)\\
&=(\|\bar{x}\|+\alpha\langle V_{p},\bar{x}\rangle)\left(\frac{\bar{x}}{\|\bar{x}\|}+\alpha V_{p}\right),\end{aligned}$$
where $\bar{x}$ is the canonical identification of $x$ (cf. p. 143 in \cite{Ziller1982}). Thus we get
$$\begin{aligned}  L_{\frac{1}{2}H_{\alpha}^{2}}\circ A^{*}(x)&=L_{\frac{1}{2}H_{\alpha}^{2}}(-x)\\
&=(\|\overline{(-x)}\|+\alpha\langle V_{-p},\overline{(-x)}\rangle)\left(\frac{\overline{(-x)}}{\|\overline{(-x)}\|}+\alpha V_{-p}\right)\\
&=(\|-\bar{x}\|+\alpha\langle -V_{p},-\bar{x}\rangle)\left(\frac{-\bar{x}}{\|-\bar{x}\|}+\alpha (-V_{p})\right)\\
&=-(\|\bar{x}\|+\alpha\langle V_{p},\bar{x}\rangle)\left(\frac{\bar{x}}{\|\bar{x}\|}+\alpha V_{p}\right)\\
&=A_{*}\circ L_{\frac{1}{2}H_{\alpha}^{2}}(x), \end{aligned}$$
where the second identity is due to $A^{*}(x)=-x\in T_{-p}^{*}S^{2n+1}$. Then (\ref{1217d}), and then
Claim 2 as well as Lemma \ref{Lm3.1} are proved. $\hfill\Box$

Here note that $M=(\mathbb{R}P^{2n+1},F)$ is orientable as well-known. Specially we have

\begin{lemma}\label{Lm3.2} Let $M=(\mathbb{R}P^{2n+1},F)$ be a Finsler manifold. Then every prime closed
geodesic on $M$ is orientable.
\end{lemma}

\Proof Note that if $\ga$ is a closed geodesic in the homotopy class $\Lm_eM$ of closed curves on $M$
represented by $e$ in $\Z_2$, it is contractible in $M$, and thus is orientable.

If $\ga$ is a closed geodesic in the homotopy class $\Lm_gM$ of closed curves on $M$ represented by $g$
in $\Z_2$. Note that $M$ can be identified with the union of the open unit ball
$$   D^{2n+1}=\{x=(x_1,\ldots,x_{2n+1})\in \R^{2n+1}\;|\; |x|< 1\}  $$
and the quotient space of the boundary $S^{2n}$ of $D^{2n+1}$ module the action of the antipodal map $A$
equipped with the Finsler metric $F$, i.e.,
\be  M = D^{2n+1}\cup (S^{2n}/A).  \lb{D.S}\ee
Thus this prime closed geodesic $\ga$ in $\Lm_gM$ can be viewed as a curve in $D^{2n+1}$ with the two end
points being the two antipodal points $p$ and $-p$ on $S^{2n}$, and $\ga$ runs from $\ga(0)=p$ to
$\ga(\tau)=-p$. Here $\ga$ can not be approximated by closed curves in the closure of $D^{2n+1}$ and is
not contractible. But the tangent space $T_pS^{2n}$ is carried to $T_{-p}S^{2n}$ by $\ga$. Therefore
$\ga$ is orientable. $\hfill\Box$.

\begin{figure}
\begin{center}
\resizebox{8cm}{8cm}{\includegraphics*[0cm,0cm][8cm,8cm]{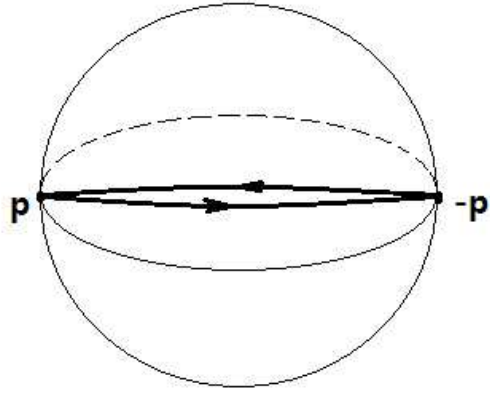}}
\vskip -5 mm
\caption{The orientability of $\ga$ and contractibility of $\ga^2$.}
\end{center}
\end{figure}
%\vspace{2mm}

\begin{remark}\label{cg2n} Note that let $\ga$ be a closed geodesic on the Finsler manifold
$M=(\mathbb{R}P^{2n},F)$ belonging to $\Lm_gM$. Then $\ga$ is not orientable.
\end{remark}

By Lemma \ref{Lm3.1}, we can endow $\mathbb{R}P^{2n+1}$ a Finsler metric induced by $N_{\alpha}$, which is still
denoted by $N_{\alpha}$ for simplicity. Therefore the natural projection
$$  \pi:(S^{2n+1},N_{\alpha})\to (\mathbb{R}P^{2n+1},N_{\alpha}),  $$
is locally isometric.

\begin{proposition}\label{cgbijiao}
For $M=(\mathbb{R}P^{2n+1},N_{\alpha})$ with $\alpha\in (0,1)\bs\Q$ and $n\geq 1$,
$$  {}^{\#}\CG(\Lm_{g}M)\leq {}^{\#}\CG(\Lm_{e}M)\leq {}^{\#}\CG(\Lm (S^{2n+1},N_{\alpha}))=2(n+1),  $$
where ${}^{\#}CG(A)$ is the number of $CG(A)$ defined in Definition \ref{def-CG(A)}, and $(S^{2n+1},N_{\alpha})$
denotes the sphere $S^{2n+1}$ endowed with the Katok metric $N_{\alpha}$.
\end{proposition}

\Proof Note that by Lemma \ref{Lm3.2}, every closed geodesic on $M$ is orientable.

For every prime closed geodesic $\ga\in \Lm_{g}M$, the second iterate $\ga^2$ of $\ga$ defined
by $\ga^{2}(t)=\ga(2t)$ is a prime closed geodesic in $\Lm_{e}M$ by the fact $g^2=e$ in $\Z_2$. This fact can be
intuitively understood as follows. As in the proof of Lemma \ref{Lm3.2}, $M$ can be identified as
in (\ref{D.S}). Thus a prime closed geodesic $\ga$ on $M$ which belongs to $\Lm_gM$ can be viewed as a
curve in $D^{2n+1}$ with the two end points being the two antipodal points $p$ and $-p$ on $S^{2n}$. Such
a closed curve in $M$ can not be approximated by closed curves in the open ball $D^{2n+1}$ and is not
contractible. Then $\ga^2$ runs from $p$ to $-p$ along $\ga$ first and then runs from $-p$ to $p$ along
$\ga^{-1}$. Now this forms a closed curve in the closure $\ol{D^{2n+1}}$ of $D^{2n+1}$, and thus we can
shrink this curve slightly at $p$ and $-p$ into $D^{2n+1}$ to form a new slightly shorter closed curve in
$D^{2n+1}$, and then continue to shrink it to the origin. We refer to Figure 3.1 for an illustration on
this argument, where for clarity, we have drawn $\ga$ and $\ga^{-1}$ slightly away from each other. From
this argument, the first inequality claimed by the proposition follows.

Let $\beta:[0,1]\to M$ be a prime closed geodesic in $\Lm_{e}M$ with $\beta(0)=\beta(1)=p$ and
$\beta^{\prime}(0)=\beta^{\prime}(1)$.  Since $(S^{2n+1}, N_{\alpha})$ is complete, we get by Hopf-Rinow theorem
there exists ${\gamma}:[0,+\infty)\to S^{2n+1}$ satisfying ${\gamma}(0)=p$ and ${\gamma}^{\prime}(0)=\beta^{\prime}(0).$
Let $\tilde{\beta}=\pi\circ{\gamma}.$ Since $\pi:(S^{2n+1},N_{\alpha})\to (\mathbb{R}P^{2n+1},N_{\alpha})$ is locally
isometric, $\tilde{\beta}$ is a geodesic in $M$ satisfying that $\tilde{\beta}(0)={\beta}(0)$ and
$\tilde{\beta}^{\prime}(0)={\beta}^{\prime}(0).$ By the uniqueness of the geodesic, $\tilde{\beta}|_{[0,1]}=\beta.$
Since the pre-image of every closed curve in $\Lm_{e}M$ is also a closed one in $S^{2n+1}$ and
$\beta^{\prime}(0)=\beta^{\prime}(1)$ implies $\gamma^{\prime}(0)=\gamma^{\prime}(1),$ $\gamma={\gamma}|_{[0,1]}$ is
a closed geodesic in $S^{2n+1}$. Note that here $\gamma$ is prime in $\Lm S^{2n+1}$ provided $\beta$ is prime in
$\Lm_e M$. Thus the second inequality follows. $\hfill\Box$

\begin{remark} (i) Here if $\ga$ is a prime closed geodesic in $(M,F)$ belonging to $\CG(\Lm_gM)$,
then $\ga^2 \in \CG(\Lm_eM)$ is also prime in $\CG(\Lm_eM)$ according to the definition \ref{def-CG(A)}. Note that $\CG(\Lm_eM)$ contains two
kinds of closed geodesics: the ones which are two iterates of $\ga$ in $\CG(\Lm_gM)$ and the ones
which are not so.

(ii) Although Proposition \ref{cgbijiao} is enough for our arguments in the sequel, we point out that the inequalities
therein  are actually equalities. Since the closed geodesics found in
\cite{Ziller1982} are great circles of $S^{2n+1}$ and the projection map $\pi$ is locally isometric, their images
under $\pi$ are exactly the two iterations of the associated non-contractible closed geodesics on $\mathbb{R}P^{2n+1}$
which implies in turn that
$${}^{\#}\CG(\Lm_{g}M)\geq{}^{\#}\CG(\Lm S^{2n+1})$$
and so they are equal.
\end{remark}

\begin{lemma}\label{Lm3.3} For $M=\mathbb{R}P^{d}$ with $d\geq2$, let $\Lambda_{g}M$ be the component of $\Lambda M$
whose elements are homotopic to $g$. Then,
$H_{*}^{S^1}(\Lambda_{g}M;\mathbb{Z}_{2})\cong H_{*}(\Lambda_{g}M/S^{1};\mathbb{Z}_{2}).$
\end{lemma}

\Proof  Let  $\Lambda[q]$ be the subspace of $\Lambda_{g}M$ with multiplicity $q$, then we have
$$\Lambda_{g}M=\bigcup_{q\in2\mathbb{N}-1}\Lambda[q].$$

Consider the projection map $\pi: \Lambda_{g}M\times_{S^{1}}ES^{1}\to \Lambda_{g}M/{S^{1}}$ defined by
$$\pi(l,e)=[l],\ \forall (l,e)\in \Lambda_{g}M\times_{S^{1}}ES^{1},$$
where $[l]$ is the equivalent class of $l$ under the natural action of $S^{1}$.

If $(l,e)\in \Lambda[q]\times_{S^1}ES^{1}$, then
$$(l,e)=\left(l\cdot \frac{k}{q},\frac{k}{q}\cdot e\right)=\left(l,\frac{k}{q}\cdot e\right),\ \forall\ 0\leq k\leq q-1.$$
Therefore, for every $[x]\in\Lambda_{g}M/{S^{1}}$ with the multiplicity of $q\in 2\mathbb{N}-1,$ we have
$$ \pi^{-1}([x])\simeq ES^{1}/\mathbb{Z}_{q}. $$
Since $q$ is odd, we obtain by Theorem III 2.4 in Bredon \cite{Bre1972} that
\begin{equation}\label{sichuan1}
H_{k}(\pi^{-1}([x]);\mathbb{Z}_{2})=H_{k}(ES^{1}/\mathbb{Z}_{q};\mathbb{Z}_{2})\cong
H_{k}(ES^{1};\mathbb{Z}_{2})^{\mathbb{Z}_{q}}=\left\{\begin{array}{ll}
\mathbb{Z}_{2},&\ k=0,\\
0,&\ \text{otherwise.}\end{array}\right.
\end{equation}
So for the map $\pi: \Lambda_{g}M\times_{S^{1}}ES^{1}\to \Lambda_{g}M/{S^{1}}$, there exists a Leray
spectral sequence $\{E^{r},d_{r}\}$ converging to $H_{*}(\Lambda_{g}M\times_{S^{1}}ES^{1};\mathbb{Z}_{2})$
with
$$ E_{s,t}^{2} =
  H_{s}\left(\Lambda_{g}M/{S^{1}};\bigcup_{[x]\in\Lambda_{g}M/{S^{1}}}H_{t}(\pi^{-1}([x]);\mathbb{Z}_{2})\right), $$
where $\bigcup_{[x]\in\Lambda_{g}M/{S^{1}}}H_{t}(\pi^{-1}([x]);\mathbb{Z}_{2})$ is the coefficient sheaf
on $\Lambda_{g}M/{S^{1}}$ with germ $H_{t}(\pi^{-1}([x]);\mathbb{Z}_{2}))$ (cf. \cite{Borel1960}, p. 35;
or \cite{Hsiang1975}, p. 37).

By (\ref{sichuan1}), the elements of $E^{2}$ are all zeroes except for the bottom line and so the spectral
sequence $\{E^{r},d_{r}\}$ collapses at $E^2$. As a result, we get that
$H_{*}^{S^1}(\Lambda_{g}M;\mathbb{Z}_{2})\cong H_{*}(\Lambda_{g}M/S^{1};\mathbb{Z}_{2})$. $\hfill\Box$

From now on, we always assume that every closed geodesic $c\in \Lambda_{g}M$ satisfies the isolated
condition (ISO) in Section 1. Note that by Lemma \ref{Lm3.2}, such a geodesic $c$ is orientable, and
its odd iterates stay still in $\Lambda_{g}M$ and even iterates stay in $\Lambda_eM$. By Theorem 1.1 of
\cite{Liu2005}, the Bott-type iteration formula Theorem 9.2.1 on p.199 and Lemma 5.3.1 on p.120 of
\cite{Lo2002} can be applied to every closed geodesic $c\in \Lm_gM$ to get
\be  i(c^{2m-1}) = \sum_{\om^{2m-1}=1}i_{\om}(c) = i(c) + 2\sum_{0<k\le m}i_{e^{k\pi/(2m-1)}}(c)
      = i(c) \quad \mod\;2,   \lb{Bott.1}\ee
for all $m\in\N$. We also mention that the iterated index formulaes under the rational coefficient field
$\mathbb{Q}$ in Section 3 of \cite{BL2010}
still hold for $\mathbb{Z}_{2}$ in this paper. Indeed, the rational coefficient field $\mathbb{Q}$ is
crucially used in the proof of Lemma 3.6 in \cite{BL2010} to ensure Theorem III 7.2 of Bredon \cite{Bre1972}
which also holds for the coefficient field $\mathbb{Z}_{2}$ since the multiplicity of every closed curve in
$\Lambda_{g}M$ is odd.

\begin{lemma}\label{Lm3.4} Let $(M,F)$ be a compact Finsler manifold with finite fundamental group and
possess only finitely many distinct prime closed geodesics, among which the homologically visible ones
are denoted by $c_i$ for $1\le i\le k$. Then we have
\be  \hat{i}(c_i) > 0, \qquad \forall\;1\le i\le k.  \lb{fcg.1}\ee
\end{lemma}

\Proof It is well known that every closed geodesic $c$ on $M$ must have mean index $\hat{i}(c)\ge 0$.

Assume by contradiction that there is a homologically visible closed geodesic $c$ on $M$ satisfying
$\hat{i}(c)=0$. Then $i(c^m)=0$ for all $m\in\N$ by Bott iteration formula and $c$ must be an absolute
minimum of $E$ in its free homotopy class, since otherwise there would exist infinitely many distinct
closed geodesics on $M$ by Theorem 3 on p.385 of \cite{BK1983}. It also follows that this homotopy
class must be non-trivial.

On the other hand, by Lemma 7.1 of \cite{Rad1992}, there exists an integer $k(c)>0$ such that
$\nu(c^{m+k(c)})=\nu(c^m)$ for all $m\in\N$. Specially we obtain $\nu(c^{mk(c)+1}) = \nu(c)$ for all
$m\in\N$ and then elements of $\ker(E''(c^{mk(c)+1}))$ are precisely $mk(c)+1$st iterates of
elements of $\ker(E''(c))$. Thus by the Gromoll-Meyer theorem in \cite{GM1969Top}, the behavior of the
restriction of $E$ to $\ker(E''(c^{mk(c)+1}))$ is the same as that of the restriction of $E$ to
$\ker(E''(c))$. Then together with the fact $i(c^m)=0$ for all $m\in\N$, we obtain that $c^{mk(c)+1}$
is a local minimum of $E$ in its free homotopy class for every $m\in\N$. Because $M$ is compact and
possessing finite fundamental group, there must exist infinitely many distinct closed geodesics on $M$
by Corollary 2 on p.386 of \cite{BK1983}. Then it yields a contradiction and proves (\ref{fcg.1}).
$\hfill\Box$.

By Lemmas \ref{Lm3.2} and \ref{Lm3.4}, we have the following result.

\begin{proposition}\label{Bettibound} For $M=\mathbb{R}P^{2n+1}$ with $n\geq1$, the $S^1$-equivariant
Betti number sequence \break $\{\bar{\beta}_k(\Lm_gM;\mathbb{Z}_{2})\}$ of $H_{*}^{S^1}(\Lambda_{g}M;\mathbb{Z}_{2})$
is bounded.
\end{proposition}

\Proof Since the $S^1$-equivariant Betti numbers $\bar{\bb}_k(\Lm_gM;\mathbb{Z}_{2})$ are topological invariants
of $M$, they are independent of the choice of the Finsler metric $F$ on it. To estimate them, it suffices to
choose a special Finsler metric $F=N_{\aa}$ for $\aa\in (0,1)\bs\Q$, i.e., the Katok metrics. Then for
$M=(\mathbb{R}P^{2n+1},N_{\alpha})$ by Proposition \ref{cgbijiao}, $\Lm_gM$ has only finitely many distinct
prime non-contractible closed geodesics, among which we denote the homologically visible ones by $c_1$, $\cdots$,
$c_r$ for some integer $r\in\N$. Note that by Lemma \ref{Lm3.4}, each $c_j$ must have positive mean index
$\hat{i}(c_j)>0$ for $1\le j\le r$. Note also that homologically invisible closed geodesics make no contributions
to any Morse type numbers. Here for $k\in\Z$, the $k$-th Morse type number $M_k(\Lm_gM)$ of $\Lm_gM$ is defined by
\be  M_k(\Lm_gM) = \sum_{j=1}^r\sum_{m\in\N}\dim\bar{C}_k(E,c_j^{2m-1}).  \lb{M.1}\ee
Now it is well known that the following inequalities hold for some constant $B>0$,
\be \bar{\beta}_k(\Lm_gM) \le M_k(\Lm_gM) \le B,  \qquad \forall\;k\in\Z, \lb{M.2}\ee
where the first inequality follows from the classical Mores theory, and the second inequality follows from
the same argument in the proof of Theorem 4 of \cite{GM1969JDG} for homologically visible prime closed geodesics
$c_1$, $\cdots$, $c_r$, and the proof is complete. $\hfill\Box$

\setcounter{figure}{0}
\setcounter{equation}{0}
\section{Possible B-V algebraic structures of
$(\mathbb{H}_{*}(L\mathbb{R}P^{2n+1};\mathbb{Z}_{2}),\bullet,\bigtriangleup)$}%Section 4

Let $M=\mathbb{R}P^{d}$ and $\pi_{1}(M)=\mathbb{Z}_{2}=\{e, g\}$ with the generator $g$ satisfying
$g^{2}=e$. Then the free loop space $LM$ decomposes into
$$  LM=L_{e}M\sqcup L_{g}M,  $$
with $L_{e}M$ and $L_{g}M$ being the connected components of $LM$ whose elements are
homotopic to $e$ and $g$ respectively. As in p.9 of \cite{CS1999}, we use the following convention:
$\mathbb{H}_{*}(LM;\mathbb{Z}_{2}) = H_{*+d}(LM;\mathbb{Z}_{2}).$

In \cite{West2005} and later in \cite{West2007}, Westerland obtained the following results.
\begin{proposition}\label{west05}
For $M=\mathbb{R}P^{2n+1}$ with $n\ge 1$,
$$  \mathbb{H}_{*}(LM; \mathbb{Z}_{2})\cong \mathbb{Z}_{2}[x,v,w]/(x^{2n+2},v^{2}-(n+1)wx^{2n}),  $$
where the topological degree of $x$, $v$ and $w$ are $-1,$ $0$ and $2n$ respectively, i.e., $|x|=-1$, $|v|=0$
and $|w|=2n.$ Moreover,
\begin{equation}
\label{homologyiso} \mathbb{H}_{*}(L_{e}M; \mathbb{Z}_{2})\cong \mathbb{H}_{*}(L_{g}M; \mathbb{Z}_{2}).
\end{equation}
\end{proposition}

\Proof For readers convenience, here we briefly indicate why (\ref{homologyiso}) should be true. From the
proof of Lemma 5.4 in \cite{West2005}, the Eilenberg-Moore spectral sequence converging to $H^{*}(LM;\mathbb{Z}_{2})$
collapses at the second page, that is, the two Eilenberg-Moore spectral sequences converging to
$H^{*}(L_{e}M;\mathbb{Z}_{2})$ and $H^{*}(L_{g}M;\mathbb{Z}_{2})$ respectively collapse at the second pages which
are the same due to (1) of Corollary 5.3 (alternatively, (1) of Lemma 5.2) therein. As a result,
we get $${H}^{*}(L_{e}M; \mathbb{Z}_{2})\cong {H}^{*}(L_{g}M; \mathbb{Z}_{2})$$ and so (\ref{homologyiso})
follows. \hfill$\Box$

Note that it is not obvious where the homology of the two components of the free loop space fit into the above
presentation of the homology of $LM.$

\begin{proposition}\label{west07}
For $M=\mathbb{R}P^{2n+1}$ with $n\ge 1$, the Poincar\'e series of $H_{*}^{S^1}(LM;\mathbb{Z}_{2})$ is
$$  \frac{1-t^{2n+2}}{(1-t^{2n})(1-t^{2})}\left(1+\frac{1+t}{1-t^{2}}\right).  $$
\end{proposition}

\Proof See Theorem 1.1 (2) in \cite{West2007}. \hfill$\Box$

We need to extract from this Poincar\'e series the part coming from the non-contractible component. For this
purpose, we first analyze the ring structure of $(\mathbb{H}_{*}(LM;\mathbb{Z}_{2}),\bullet)$.

\begin{lemma}\label{algebraic1}
With the same notations in Proposition \ref{west05},  there exist
\bea
&& \tilde{x}=x+a_{1}xv+a_{2}x^{2n+1}w+a_{3}x^{2n+1}vw\in\mathbb{H}_{-1}(M;\mathbb{Z}_{2})\subseteq\mathbb{H}_{-1}(L_{e}M;\mathbb{Z}_{2}), \nn\\
&& \tilde{v}=v+b_{1}+b_{2}\td{x}^{2n}w+b_{3}\td{x}^{2n}vw\in\mathbb{H}_{0}(L_{g}M;\mathbb{Z}_{2}),  \nn\eea
and
$$\tilde{w}=c_{0}w+c_{1}\td{v}w+c_{2}\td{x}^{2n}w^{2}+c_{3}\td{x}^{2n}\td{v}w^{2}
            \in  \mathbb{H}_{2n}(L_{e}M;\mathbb{Z}_{2})\cup\mathbb{H}_{2n}(L_{g}M;\mathbb{Z}_{2}), $$
where $a_{j}, b_{j},c_{j}\in\mathbb{Z}_{2}$ are some coefficients found in the subsequent proofs of this lemma, such that
\begin{equation}
\label{summarycong}
\mathbb{H}_{*}(LM; \mathbb{Z}_{2})
\cong
 \mathbb{Z}_{2}[\tilde{x},\tilde{v},\tilde{w}]/(\tilde{x}^{2n+2},\tilde{v}^{2}+b_{1}
             -(n+1)\sigma(b_{1},c_{1})\tilde{x}^{2n}\tilde{v}^{b_{1}c_{1}}\tilde{w}),
 \end{equation}
where
$$\sigma(b_{1},c_{1})=
 \left\{
 \begin{array}{ll}
 0,&\ \text{if}\ b_{1}=0\ \text{and}\ c_{1}=1;\\
 1,&\ \text{otherwise.}
 \end{array}
 \right. $$
\end{lemma}

\Proof We construct $\tilde{x}$, $\tilde{v}$ and $\tilde{w}$ one by one as follows.

{\bf Step 1: Construction of $\tilde{x}$.}

Consider first the composed maps
$$  M\stackrel{i}{\longrightarrow}L_{e}M\stackrel{\ev}{\longrightarrow}M,  $$
where $i$ is the inclusion map which maps each point in $M$ to the corresponding constant loop in $L_eM$ and
$\ev(\gamma)=\gamma(0)$ is the evaluation map at the initial point of the curve. Since $\ev\circ i = \id$, the
map
$$  i_{*}:\mathbb{H}_{*}(M;\mathbb{Z}_{2})\to \mathbb{H}_{*}(L_{e}M;\mathbb{Z}_{2}),  $$
is an embedding.

Let $\tilde{x}\in\mathbb{H}_{-1}(M;\mathbb{Z}_{2})$ be the generator of $\mathbb{H}_{*}(M;\mathbb{Z}_{2})$. It follows
by Proposition \ref{west05} that there exist $a_{0},a_{1},a_{2}, a_{3}\in \mathbb{Z}_{2}$ such that
$$  \tilde{x} =a_{0} x+a_{1}xv+a_{2}x^{2n+1}w+a_{3} x^{2n+1}vw.  $$
We claim $a_{0}=1$, that is,
 \begin{equation}
 \label{2014714a}
 \tilde{x} = x+a_{1}xv+a_{2}x^{2n+1}w+a_{3} x^{2n+1}vw.
 \end{equation}
If not, then $\tilde{x} =a_{1}xv+a_{2}x^{2n+1}w+a_{3} x^{2n+1}vw.$ By Proposition \ref{west05}, $x^{2n+2}=0$
and $v^{2}=(n+1)x^{2n}w$ which imply
$$  \tilde{x}^{2}=x^{2}(a_{1}v+a_{2}x^{2n}w+a_{3} x^{2n}vw)^{2}=0,  $$
but this contradicts to $\mathbb{H}_{*}(M;\mathbb{Z}_{2})\cong \mathbb{Z}_{2}[\tilde{x}]/(\tilde{x}^{2n+2})$.

Further direct computations show that for $k\geq2,$
$$\begin{aligned}
\td{x}^k
&= (x + a_{1}xv + a_{2}x^{2n+1}w + a_{3} x^{2n+1}vw)^{k}  \nn\\
&=x^k(1 + a_1v + a_{2}x^{2n}w + a_{3} x^{2n}vw)^{k}  \nn\\
&= x^k(1 + a_1v)^{k}  \nn\\
&= x^{k}(1 + ka_1v),\end{aligned}  $$
where the third identity is due to $x^{2n+k}=0$ and the forth one follows by $x^{k}v^{2}=(n+1)x^{2n+k}w=0$.
Therefore we obtain that for $k\geq2$,
\be  \td{x}^k =
\left\{\begin{array}{cc}
      x^k, & {\rm if}\;k\;{\rm is\;even}\;{\rm or}\;a_1=0, \\
      x^{k}(1 + v), & {\rm if}\;k\;{\rm is\;odd}\,{\rm and}\;a_1=1,\\
   \end{array}\right.  \lb{2014714c}\ee
which implies that
\begin{equation}
\label{2014714d}
\tilde{x}^{2n+2}=0\ \text{and}\  v^{2}-(n+1)w\tilde{x}^{2n}=0.
\end{equation}

We come back to consider the case of $k=1$. Multiplying both sides of (\ref{2014714a}) by $v$, we obtain
$$\begin{aligned}
\tilde{x}v&=xv+a_{1}xv^{2}+a_{2}x^{2n+1}vw\\
&=xv+(n+1)a_{1}x^{2n+1}w+a_{2}x^{2n+1}vw,
\end{aligned}$$
which is equivalent to
\begin{equation}\label{20150201a}
xv=\tilde{x}v+(n+1)a_{1}x^{2n+1}w+a_{2}x^{2n+1}vw.
\end{equation}
By direct computations we then get
\begin{equation}
\label{20150201c}
\begin{aligned}
x&=\tilde{x}+\td{a}_{1}\tilde{x}v+\td{a}_{2}\tilde{x}^{2n+1}w+\td{a}_{3}\tilde{x}^{2n+1}vw,
\end{aligned}
\end{equation}
where $\td{a}_{1}$, $\td{a}_{2}$ and $\td{a}_{3} \in \Z_{2}$ can be determined from ${a}_{1}$, ${a}_{2}$ and ${a}_{3}$
via (\ref{2014714a}), (\ref{2014714c}) and (\ref{20150201a}). That is, $x$ can be represented in terms of $\td{x}$, $v$
and $w$.

Due to  (\ref{20150201c}) and (\ref{2014714d}), the Chas-Sullivan ring  $\mathbb{Z}_{2}[x,v,w]/(x^{2n+2}, v^{2}-(n+1)wx^{2n})$
can also be represented in terms of $\td{x}$, $v$, $w$, and elements in the set
$$  K=\{\td{x}^{a}v^{b}w^{c}\mid0\leq a\leq 2n+1,\ 0\leq b\leq 1\ \text{and}\ 0\leq c<+\infty\},  $$
where elements in $K$ from degree $-(2n+1)$ to $0$ can be listed as follows:
\begin{equation}
\label{degreelist}\left\{
\begin{array}{ll}
\deg=\ 0;&\ 1,\ v,\ \td{x}^{2n}w,\ \td{x}^{2n}vw;\\
\deg=-1;&\ \td{x},\ \td{x}v,\ \td{x}^{2n+1}w,\ \td{x}^{2n+1}vw;\\
\deg=-2;&\ \td{x}^{2},\ \td{x}^{2}v;\\
\deg=-3;&\ \td{x}^{3},\ \td{x}^{3}v;\\
\qquad \vdots&\qquad \ \vdots\\
\deg=-2n+1;&\ \td{x}^{2n-1},\ \td{x}^{2n-1}v;\\
\deg=-2n;&\ \td{x}^{2n},\ \td{x}^{2n}v;\\
\deg=-(2n+1);&\ \td{x}^{2n+1},\ \td{x}^{2n+1}v;\\
\end{array}\right.
\end{equation}
while the elements in $K$ from degree $1$ to $2n$ are the corresponding ones from degree $-2n+1$ to $0$
in (\ref{degreelist}) multiplying by $w$; the elements in $K$ with higher degrees can be
obtained by the similar way via multiplying $w^k$ for integers $k>1$. This representation is an isomorphism by the
following claim.

{\bf Claim:} {\it The Chas-Sullivan ring
$\mathbb{H}_{*}(LM; \mathbb{Z}_{2}) \cong \mathbb{Z}_{2}[x,v,w]/(x^{2n+2},  v^{2}-(n+1)wx^{2n})$ can be
represented isomorphically in terms of $\tilde{x}$, $v$, and $w$ as
\begin{equation}
\label{cong1}   \mathbb{H}_{*}(LM; \mathbb{Z}_{2}) \cong \mathbb{Z}_{2}[\tilde{x},v,w]/(\tilde{x}^{2n+2}, v^{2}
       -(n+1)\tilde{x}^{2n}w). \end{equation}}

To prove this claim, it suffices to show that the elements in $K$ with the same degree are linearly independent.

Notice also that  $yw=0$ implies $y=0$ for any $y\in\mathbb{Z}_{2}[x,v,w]/(x^{2n+2},  v^{2}-(n+1)wx^{2n})$. Indeed,
assume without loss of generality $y=\sum_{1\leq i\leq i_{0}}k_{i}{x}^{a_{i}}v^{b_{i}}w^{c_{i}}$ with $i_{0}$ the
dimension of the linear space with degree $|y|$ , $k_{i}\in\Z_{2}$ and $-a_{i}+2n c_{i}=|y|$.  Then, we get
$$  \sum_{1\leq i\leq i_{0}}k_{i}{x}^{a_{i}}v^{b_{i}}w^{c_{i}+1}=yw=0.  $$
But the elements ${x}^{a_{i}}v^{b_{i}}w^{c_{i}+1}$ with ${1\leq i\leq i_{0}}$ are linearly independent due to the
ring structure itself, which implies $k_{i}=0$ for every $1\leq i\leq i_{0}$ and thus $y=0$.

As a result, we need only to prove that the elements with the same degree listed in (\ref{degreelist}) are
linearly independent.

$1^{\circ}$ For $\deg=-1$, assume that there is
$  (l_{1},l_{2},l_{3},l_{4})\in\mathbb{Z}_{2}\times\mathbb{Z}_{2}\times\mathbb{Z}_{2}\times\mathbb{Z}_{2} $
such that
\begin{equation}
\label{2014726a}
l_{1}\tilde{x}+l_{2}\tilde{x}v+l_{3}\tilde{x}^{2n+1}w+l_{4} \tilde{x}^{2n+1}vw=0.
\end{equation}
Multiplying both sides of (\ref{2014726a}) with $x$, we get by (\ref{2014714a}) and (\ref{2014714c}) that
$$0=l_{1}x(x+a_{1}xv)+l_{2}x^{2}v=l_{1}x^{2}+(a_{1}l_{1}+l_{2})x^{2}v,$$
which implies $l_{1}=l_{2}=0.$

Therefore, (\ref{2014726a}) can be rewritten as
\begin{equation}
\label{2014726b}
l_{3}\tilde{x}^{2n+1}w+l_{4} \tilde{x}^{2n+1}vw=0.
\end{equation}
Again by (\ref{2014714c}), we have
$$\begin{aligned}l_{3}\tilde{x}^{2n+1}w+l_{4} \tilde{x}^{2n+1}vw&=\left\{\begin{array}{ll}
     l_{3}x^{2n+1}w+l_{4} x^{2n+1}vw, & {\rm if}\;a_1=0, \\
     l_{3}x^{2n+1}(1+v)w+l_{4} x^{2n+1}(1+v)vw, & {\rm if}\;a_1=1,\\
   \end{array}\right.\\
   &=\left\{\begin{array}{ll}
     l_{3}x^{2n+1}w+l_{4} x^{2n+1}vw, & {\rm if}\;a_1=0, \\
     l_{3}x^{2n+1}w+(l_{3}+l_{4}) x^{2n+1}vw, & {\rm if}\;a_1=1.\\
   \end{array}\right.\\
   \end{aligned}$$
which implies $l_{3}=l_{4}=0$ too.

$2^{\circ}$ For $\deg=-(2k+1)$ with $1\leq k\leq n$, the proof is similar as above. Assume that there is
$(l_{1},l_{2})\in\mathbb{Z}_{2}\times\mathbb{Z}_{2}$
such that
\begin{equation}
\label{2014823a}
l_{1}\tilde{x}^{2k+1}+l_{2}\tilde{x}^{2k+1}v=0.
\end{equation}
Then by (\ref{2014714c}) we get
$$
\begin{aligned}
0&=\left\{\begin{array}{cc}
      l_{1}{x}^{2k+1}+l_{2}{x}^{2k+1}v, & {\rm if}\;a_1=0, \\
      l_{1}x^{2k+1}(1+ v)+l_{2}{x}^{2k+1}(1+v)v, & {\rm if}\;a_1=1,\\
   \end{array}\right. \\
    &=\left\{\begin{array}{cc}
      l_{1}{x}^{2k+1}+l_{2}{x}^{2k+1}v, & {\rm if}\;a_1=0, \\
      l_{1}x^{2k+1}+(l_{1}+l_{2}){x}^{2k+1}v, & {\rm if}\;a_1=1,\\
   \end{array}\right.
   \end{aligned}
$$
which implies $l_{1}=l_{2}=0.$

$3^{\circ}$ For $\deg=-2k$ with $0\leq k\leq n$, the proof follows from the linear independence of $x^{2k}$ and $x^{2k}v$
when $k>0$, and that of $1$, $v$, $x^{2n}w$, and $x^{2n}vw$ when $k=0$, due to $\td{x}^{2k}={x}^{2k}$ by (\ref{2014714c}).
Then the claim is proved.

{\bf Step 2: Construction of $\tilde{v}$.}

Again by Proposition \ref{west05},  each of $\mathbb{H}_{-(2n+1)}(L_{e}M;\mathbb{Z}_{2})$ and
$\mathbb{H}_{-(2n+1)}(L_{g}M;\mathbb{Z}_{2})$ has a generator in $\{\tilde{x}^{2n+1},\tilde{x}^{2n+1}v, \tilde{x}^{2n+1}(1+v)\}$.
Since
$$  \tilde{x}^{2n+1}\in\mathbb{H}_{-(2n+1)}(M;\mathbb{Z}_{2})\subseteq\mathbb{H}_{-(2n+1)}(L_{e}M;\mathbb{Z}_{2}),  $$
either $\tilde{x}^{2n+1}v$ or $\tilde{x}^{2n+1}(1+v)$ belongs to $\mathbb{H}_{-(2n+1)}(L_{g}M;\mathbb{Z}_{2})$. According to
the definition of the loop product (cf. pp. 6-7, in \cite{CS1999}), we have
\bea
\mathbb{H}_{*}(L_{e}M;\mathbb{Z}_{2})\bullet\mathbb{H}_{*}(L_{e}M;\mathbb{Z}_{2}) &\subset& \mathbb{H}_{*}(L_{e}M;\mathbb{Z}_{2}), \lb{Map1}\\
\mathbb{H}_{*}(L_{e}M;\mathbb{Z}_{2})\bullet\mathbb{H}_{*}(L_{g}M;\mathbb{Z}_{2}) &\subset& \mathbb{H}_{*}(L_{g}M;\mathbb{Z}_{2}), \lb{Map2}\\
\mathbb{H}_{*}(L_{g}M;\mathbb{Z}_{2})\bullet\mathbb{H}_{*}(L_{g}M;\mathbb{Z}_{2}) &\subset& \mathbb{H}_{*}(L_{e}M;\mathbb{Z}_{2}). \lb{Map3}\eea
So there exists
$$ \tilde{v}\in\mathbb{H}_{0}(L_{g}M;\mathbb{Z}_{2})\subset\mathbb{H}_{0}(LM;\mathbb{Z}_{2})=\Span\{1,v,\tilde{x}^{2n}w,\tilde{x}^{2n}vw\} $$
such that
$$ 0\neq\tilde{x}^{2n+1}\tilde{v}\in\mathbb{H}_{-(2n+1)}(L_{g}M;\mathbb{Z}_{2}), $$
and we assume
$$ \tilde{v}=b_{0}v+b_{1}+b_{2}\tilde{x}^{2n}w+b_{3}\tilde{x}^{2n}vw, $$
with $b_{0},b_{1},b_{2},b_{3}\in\mathbb{Z}_{2}$.

Now we claim $b_{0}=1$, that is
\begin{equation}
\label{2014715a}
\tilde{v}=v+b_{1}+b_{2}\tilde{x}^{2n}w+b_{3}\tilde{x}^{2n}vw.
\end{equation}
If not, then we get
$$\tilde{x}^{2n+1}\tilde{v}=\tilde{x}^{2n+1}(b_{1}+b_{2}\tilde{x}^{2n}w+b_{3}\tilde{x}^{2n}vw)
   =b_{1}\tilde{x}^{2n+1}\in\mathbb{H}_{-(2n+1)}(L_{e}M;\mathbb{Z}_{2}),$$
a contradiction.

It then follows from (\ref{2014715a}) that
\begin{equation}
\label{2014715b}
\tilde{v}^{2}=\left\{
\begin{array}{ll}
v^{2},&\ \text{if}\ b_{1}=0,\\
v^{2}+1,&\ \text{if}\ b_{1}=1,\\
\end{array}
\right.
\end{equation}
and
\begin{equation}
\label{2014715c}
\begin{aligned}
v&=\tilde{v}+b_{1}+b_{2}\tilde{x}^{2n}w+b_{3}\tilde{x}^{2n}vw\\
&=\tilde{v}+b_{1}+b_{2}\tilde{x}^{2n}w+b_{3}\tilde{x}^{2n}(\tilde{v}+b_{1}+b_{2}\tilde{x}^{2n}w+b_{3}\tilde{x}^{2n})w\\
&=\tilde{v}+b_{1}+(b_{2}+b_{1}b_{3})\tilde{x}^{2n}w+b_{3}\tilde{x}^{2n}\tilde{v}w.\\
\end{aligned}
\end{equation}
That is, $v$ can be represented in terms of $\td{x}$, $\td{v}$ and $w$.

Now by (\ref{2014715a}), (\ref{2014715b}) and (\ref{2014715c}), similar arguments as those in Step 1 then yield
\begin{equation}
  \label{cong2}
  \begin{aligned}
\mathbb{Z}_{2}[\tilde{x},v,w]/(\tilde{x}^{2n+2},{v}^{2}-(n+1)w\tilde{x}^{2n})
\cong
 \mathbb{Z}_{2}[\tilde{x},\tilde{v},w]/(\tilde{x}^{2n+2},\tilde{v}^{2}+b_{1}-(n+1)w\tilde{x}^{2n}).
\end{aligned}
\end{equation}

{\bf Step 3: Construction of $\tilde{w}$.}

We consider
$$ w \in \mathbb{H}_{2n}(LM;\mathbb{Z}_{2})=\Span\{w,\tilde{v}w,\tilde{x}^{2n}w^{2},\tilde{x}^{2n}\tilde{v}w^{2}\}. $$

If
\begin{equation}\label{201484a}
{w}\in \mathbb{H}_{2n}(L_{e}M;\mathbb{Z}_{2})\cup\mathbb{H}_{2n}(L_{g}M;\mathbb{Z}_{2}),
\end{equation}
we need do nothing but define $\td{w}=w$.

While if
$$  {w}\in\mathbb{H}_{2n}(L_{e}M;\mathbb{Z}_{2})+\mathbb{H}_{2n}(L_{g}M;\mathbb{Z}_{2}),  $$
then due to the linear independence of $w$ with the other three generators, there exist
$c_{1},c_{2},c_{3 }\in\mathbb{Z}_{2}$ so that $w$ can be decomposed into two parts as follows
$$ w = (w+c_{1}\tilde{v}w+c_{2}\tilde{x}^{2n}w^{2}+c_{3}\tilde{x}^{2n}\tilde{v}w^2)+(c_{1}\tilde{v}w^{2}
         +c_{2}\tilde{x}^{2n}w^{2}+c_{3}\tilde{x}^{2n}\tilde{v}w^{2}),  $$
with one part in $\mathbb{H}_{2n}(L_{e}M;\mathbb{Z}_{2})$ and the other in $\mathbb{H}_{2n}(L_{g}M;\mathbb{Z}_{2})$.

We continue the discussion in three cases according to the possible values of $c_1$, $c_2$ and $c_3$.

(1) If $c_{1}=0$, we define
\begin{equation}
\label{2014716a}
\begin{aligned}
\tilde{w}&=w+c_{1}\tilde{v}w+c_{2}\tilde{x}^{2n}w^{2}+c_{3}\tilde{x}^{2n}\tilde{v}w^{2}\\
&=w+c_{2}\tilde{x}^{2n}w^{2}+c_{3}\tilde{x}^{2n}\tilde{v}w^{2}.\\
\end{aligned}
\end{equation}
Then, multiplying by $\tilde{x}^{2n}$ both sides of  (\ref{2014716a}) we get
\begin{equation}
\label{2014716b}
\tilde{x}^{2n}\tilde{w}=\tilde{x}^{2n}{w}
\end{equation}
 and so
\begin{equation}
\label{2014716c}
\begin{aligned}
w&=\tilde{w}+c_{2}\tilde{x}^{2n}w^{2}+c_{3}\tilde{x}^{2n}\tilde{v}w^{2}\\
&=\tilde{w}+c_{2}\tilde{x}^{2n}\tilde{w}w+c_{3}\tilde{x}^{2n}\tilde{v}\tilde{w}w\\
&=\tilde{w}+c_{2}\tilde{x}^{2n}\tilde{w}^{2}+c_{3}\tilde{x}^{2n}\tilde{v}\tilde{w}^{2}\\
\end{aligned}
\end{equation}
That is, $w$ can be represented in terms of $\td{x}$, $\td{v}$ and $\td{w}$.

By (\ref{2014716a}), (\ref{2014716b}) and (\ref{2014716c}), similar arguments as those in Step 1 then yield
\begin{equation}
\label{cong3a}
\begin{aligned}
 \mathbb{Z}_{2}[\tilde{x},\tilde{v},w]/(\tilde{x}^{2n+2},\tilde{v}^{2}+b_{1}-(n+1)w\tilde{x}^{2n})
  \cong &
 \mathbb{Z}_{2}[\tilde{x},\tilde{v},\tilde{w}]/(\tilde{x}^{2n+2},\tilde{v}^{2}+b_{1}-(n+1)\tilde{w}\tilde{x}^{2n}).
 \end{aligned}
\end{equation}

(2) If $c_{1}=1$ and $b_{1}=0$, we also define
\begin{equation}
\label{2014716aaa}
\begin{aligned}
\tilde{w}&=w+c_{1}\tilde{v}w+c_{2}\tilde{x}^{2n}w^{2}+c_{3}\tilde{x}^{2n}\tilde{v}w^{2}\\
&=w(1+\tilde{v})+c_{2}\tilde{x}^{2n}w^{2}+c_{3}\tilde{x}^{2n}\tilde{v}w^{2}.\\
\end{aligned}
\end{equation}

Then, multiplying  both sides of  (\ref{2014716aaa}) by $\tilde{x}^{2n}(1+\tilde{v})$ and $(1+\tilde{v})$ respectively  we get
\begin{equation}
\label{2014716b2}
\tilde{x}^{2n}\tilde{w}(1+\tilde{v})=\tilde{x}^{2n}{w}(1+\tilde{v})^{2}=\tilde{x}^{2n}{w},
\end{equation}
and
\begin{equation}
\label{2014717a}
\begin{aligned}
\tilde{w}(1+\tilde{v})=w(1+\tilde{v}^{2})+c_{2}\tilde{x}^{2n}(1+\tilde{v}){w}^{2}+c_{3}\tilde{x}^{2n}\tilde{v}{w}^{2},
\end{aligned}
\end{equation}
where we have used the identities $$(1+\tilde{v})^{2}=1+2\tilde{v}+\tilde{v}^{2}=1+\tilde{v}^{2}=1+{v}^{2}=1+(n+1)\tilde{x}^{2n}w.$$

Notice also that $n$ must be odd in this case. Indeed, if $n$ is even,
$\tilde{x}^{2n}\tilde{w}(1+\tilde{v})=\tilde{x}^{2n}{w}=\tilde{v}^{2}\in\mathbb{H}_{0}(L_{e}M;\mathbb{Z}_{2})$
by (\ref{2014716b2}), but both $\tilde{x}^{2n}\tilde{w}$ and $\tilde{x}^{2n}\tilde{w}\tilde{v}$ are nonzero
and belong to different parts of $\mathbb{H}_{0}(L_{e}M;\mathbb{Z}_{2})\cup \mathbb{H}_{0}(L_{g}M;\mathbb{Z}_{2})$, a
contradiction. As a result, we have
\begin{equation}
\label{20150205}
\td{v}^{2}=(n+1)\tilde{x}^{2n}w=0.
\end{equation}

It then follows by (\ref{2014716b2}), (\ref{2014717a}) and (\ref{20150205}) that
\begin{equation}
\label{2014716c2}
\begin{aligned}
w&=\tilde{w}(1+\tilde{v})+w\tilde{v}^{2}+c_{2}\tilde{x}^{2n}(1+\tilde{v})w^{2}+c_{3}\tilde{x}^{2n}\tilde{v}w^{2}\\
&=\tilde{w}(1+\tilde{v})+c_{2}\tilde{x}^{2n}(1+\tilde{v})\tilde{w}^{2}+c_{3}\tilde{x}^{2n}\tilde{v}\tilde{w}^{2}.\\
\end{aligned}
\end{equation}
That is, $w$ can be represented in terms of $\td{x}$, $\td{v}$ and $\td{w}$.

By (\ref{2014716aaa}), (\ref{20150205}) and (\ref{2014716c2}), similar arguments as those in Step 1 then yield
\begin{equation}
\label{cong3b}
\begin{aligned}
&
 \mathbb{Z}_{2}[\tilde{x},\tilde{v},w]/(\tilde{x}^{2n+2},\tilde{v}^{2}) \cong
 \mathbb{Z}_{2}[\tilde{x},\tilde{v},\tilde{w}]/(\tilde{x}^{2n+2},\tilde{v}^{2}).\\
 \end{aligned}
\end{equation}

(3) If $c_{1}=1$ and $b_{1}=1$, we define
\begin{equation}
\label{2014716aa}
\begin{aligned}
\tilde{w}&=c_{1}\tilde{v}w+c_{2}\tilde{x}^{2n}w^{2}+c_{3}\tilde{x}^{2n}\tilde{v}w^{2}\\
&=\tilde{v}w+c_{2}\tilde{x}^{2n}w^{2}+c_{3}\tilde{x}^{2n}\tilde{v}w^{2}\\.
\end{aligned}
\end{equation}

Multiplying  both sides of  (\ref{2014716aa}) by $\tilde{x}^{2n}\tilde{v}$ and $\tilde{v}$, we get
\begin{equation}
\label{2014717b2}
\tilde{x}^{2n}\tilde{w}\tilde{v}=\tilde{x}^{2n}{w}\tilde{v}^{2}=\tilde{x}^{2n}{w},
\end{equation}
and
\begin{equation}
\label{2014717b3}
\begin{aligned}
\tilde{w}\tilde{v}
&=\tilde{v}^{2}w+c_{2}\tilde{x}^{2n}\tilde{v}w^{2}+c_{3}\tilde{x}^{2n}\tilde{v}^{2}w^{2}\\
&=(\tilde{x}^{2n}w+1)w+c_{2}\tilde{x}^{2n}\tilde{v}w^{2}+c_{3}\tilde{x}^{2n}\tilde{v}^{2}w^{2}\\
&=w+\tilde{x}^{2n}w^{2}+c_{2}\tilde{x}^{2n}\tilde{v}w^{2}+c_{3}\tilde{x}^{2n}w^{2}\\
\end{aligned}
\end{equation}
By (\ref{2014717b2}) and (\ref{2014717b3}) we have
\begin{equation}
\label{2014717c2}
\begin{aligned}
w&=\tilde{w}\tilde{v}+c_{2}\tilde{x}^{2n}\tilde{v}w^{2}+(1+c_{3})\tilde{x}^{2n}w^{2}\\
&=\tilde{w}\tilde{v}+c_{2}\tilde{x}^{2n}\tilde{v}^{3}\tilde{w}^{2}+(1+c_{3})\tilde{x}^{2n}\tilde{v}^{2}\tilde{w}^{2}\\
&=\tilde{w}\tilde{v}+c_{2}\tilde{x}^{2n}\tilde{v}\tilde{w}^{2}+(1+c_{3})\tilde{x}^{2n}\tilde{w}^{2}.\\
\end{aligned}
\end{equation}
That is, $w$ can be represented in terms of $\td{x}$, $\td{v}$ and $\td{w}$.

Then by (\ref{2014716aa}), (\ref{2014717b2}) and (\ref{2014717c2}), similar arguments as those in Step 1 yield
\begin{equation}
\label{cong3c}
\begin{aligned}
 \mathbb{Z}_{2}[\tilde{x},\tilde{v},w]/(\tilde{x}^{2n+2},\tilde{v}^{2}+1-(n+1)w\tilde{x}^{2n})
 \cong
  \mathbb{Z}_{2}[\tilde{x},\tilde{v},\tilde{w}]/(\tilde{x}^{2n+2},\tilde{v}^{2}+1-(n+1)\tilde{x}^{2n}\tilde{v}\tilde{w}).
 \end{aligned}
\end{equation}

Finally, summarizing (\ref{cong1}), (\ref{cong2}), (\ref{cong3a}), (\ref{cong3b}) and (\ref{cong3c}) together,
we obtain (\ref{summarycong}) and complete the proof of the lemma. \hfill$\Box$

\begin{remark}\label{1219a}
In the sequel, we are not concerned with the precise relationship of  $\td{x}$, $\td{v}$ and $\td{w}$ in Lemma
\ref{algebraic1}. Actually whatever it is, $\mathbb{H}_{*}(LM;\mathbb{Z}_{2})$ is generated by
$$ \left\{\td{x}^{l}\td{v}^{m}\td{w}^{n}\mid 0\leq l\leq 2n+1, 0\leq m\leq 1\  {\rm and}\ 0\leq n<+\infty\right\}, $$
and we only use $\td{x}$, $\td{v}$ and $\td{w}$ to separate $\mathbb{H}_{*}(L_{e}M;\mathbb{Z}_{2})$ and
$\mathbb{H}_{*}(L_{g}M;\mathbb{Z}_{2})$ from $\mathbb{H}_{*}(LM;\mathbb{Z}_{2})$ in Theorem \ref{ThmBVstr}.
For notational simplicity, we still use $x$, $v$ and $w$ later instead of $\tilde{x}$, $\tilde{v}$ and $\tilde{w}$
respectively.
\end{remark}

We need the following lemma for the purpose of computations.

\begin{lemma}\label{fanfujisuan}
For every $x,y\in\mathbb{H}_{*}(LM;\mathbb{Z}_{2})$ and $k\in\mathbb{N}$,

$(i)$ $\bigtriangleup(x^{2})=0,$ $\{x,y^{2}\}=0.$

$(ii)$ $\bigtriangleup(x^{2k})=0,$ $\{x,y^{2k}\}=0.$

$(iii)$ $\bigtriangleup(xy^{2k})=\bigtriangleup(x)y^{2k}$  and $\{x,y^{2k+1}\}=\{x,y\}y^{2k}.$
\end{lemma}

\Proof (i) According to the definition 4.1 in \cite{CS1999} (p. 12) and observing that the coefficient
field is $\mathbb{Z}_{2}$, we have $\{x,x\}=0$. By the B-V formulae,
$\bigtriangleup(x^{2})=x\bullet\bigtriangleup(x)+\bigtriangleup(x)\bullet x+\{x,x\}=0.$ Similarly,
$$  \{x,y^{2}\}=\{x,y\}\bullet y+y\bullet\{x,y\}=0.  $$

(ii) $\bigtriangleup(x^{2k})=\bigtriangleup((x^{k})^{2})=0$, $\{x,y^{2k}\}=\{x,(y^{k})^{2}\}=0.$

(iii) By (ii),
$$ \bigtriangleup(xy^{2k}) = \bigtriangleup(x)\bullet y^{2k}+x\bullet\bigtriangleup(y^{2k})+\{x,y^{2k}\}
                = \bigtriangleup(x)y^{2k}.$$
Similarly,
$$  \{x,y^{2k+1}\}=\{x,y\bullet y^{2k}\}=\{x,y\}\bullet y^{2k}+y\bullet\{x,y^{2k}\}=\{x,y\}y^{2k}. $$
$\hfill\Box$

In the sequel, we will focus on the computation of the third pages of  $\LS(L_{e}M)$ and $\LS(L_{g}M)$ whose
second pages are $H_{*}(BS^{1})\otimes \mathbb{H}_{*}(L_{e}M)$ and $H_{*}(BS^{1})\otimes \mathbb{H}_{*}(L_{g}M)$,
respectively. As usual, we introduce the following definition for convenience of writing.
\begin{definition}
For a subset $S$ of $H_{*}(BS^{1})\otimes \mathbb{H}_{*}(L_{e}M)$ (resp. $H_{*}(BS^{1})\otimes \mathbb{H}_{*}(L_{e}M)$),
the set generated by $S$ consists of linear combinations of elements of the same topological degree in $S$
with coefficients in $\Z_2$.
\end{definition}

\begin{lemma}\label{cruciallemma}
For $M=\mathbb{R}P^{2n+1}$ with $n\geq1$, $\bigtriangleup\equiv0$ on $\mathbb{H}_{*}(L_{e}M;\mathbb{Z}_{2})$.
\end{lemma}

\Proof Assume by contradiction that there exists $y\in\mathbb{H}_{*}(L_{e}M;\mathbb{Z}_{2})$ such that
$\bigtriangleup(y)\neq0.$ Then for every $k\in\mathbb{N}\cup\{0\}$, $yw^{2k}\in\mathbb{H}_{*}(L_{e}M;\mathbb{Z}_{2})$
and by Lemma \ref{fanfujisuan} $\bigtriangleup(yw^{2k})=\bigtriangleup(y) w^{2k}\neq0$. According to Lemma \ref{xym1},
for the Leray-Serre spectral sequence $\LS(L_eM)$ we have
$$  d_{2}(u^{p}\otimes yw^{2k})=u^{p-1}\otimes \bigtriangleup(yw^{2k})\neq0,\ \forall\ p\in\mathbb{N}.  $$
Thus at least the elements generated by the set
$$  \{u^{p}\otimes yw^{2k}\mid k\in\mathbb{N}\cup\{0\},\ p\in\mathbb{N}\}  $$
are killed when the spectral sequence passes to the third page from the second one, since they are not in the
kernel of $d_{2}$. As a result, the Poincar\'e series $P_{III}^{S^{1}}(L_{e}M;\mathbb{Z}_{2})(t)$ associated to
the third page of $\LS(L_{e}M)$ satisfies
\bea P_{III}^{S^{1}}(L_{e}M;\mathbb{Z}_{2})(t)
&\le& P_{II}^{S^{1}}(L_{e}M;\mathbb{Z}_{2})(t)-t^{2n+1}t^{|y|}\sum_{p=1}^{+\infty}\sum_{k=0}^{+\infty}t^{2p}t^{4kn} \nn\\
&=& P_{II}^{S^{1}}(L_{e}M;\mathbb{Z}_{2})(t) - \frac{t^{2n+3+|y|}}{(1-t^{2})(1-t^{4n})}. \nn\eea
Therefore, we get
\bea P^{S^{1}}(L_{g}M;\mathbb{Z}_{2})(t)
&=& P^{S^{1}}(LM;\mathbb{Z}_{2})(t) - P^{S^{1}}(L_{e}M;\mathbb{Z}_{2})(t)   \nn\\
&\ge& P^{S^{1}}(LM;\mathbb{Z}_{2})(t) - P_{III}^{S^{1}}(L_{e}M;\mathbb{Z}_{2})(t)  \nn\\
&\ge& P^{S^{1}}(LM;\mathbb{Z}_{2})(t) - P_{II}^{S^{1}}(L_{e}M;\mathbb{Z}_{2})(t)
             + \frac{t^{2n+3+|y|}}{(1-t^{2})(1-t^{4n})} \nn\\
&=& \frac{1-t^{2n+2}}{(1-t^{2n})(1-t^{2})} + \frac{t^{2n+3+|y|}}{(1-t^{2})(1-t^{4n})} \nn\\
&>& \frac{t^{2n+3+|y|}}{(1-t^{2})(1-t^{4n})},  \label{201413a}\eea
where to get the second equality we have used direct computations for $P_{II}^{S^{1}}(L_{e}M;\mathbb{Z}_{2})(t)$
and Proposition \ref{west07}. Since the right hand side of (\ref{201413a}) contains the multiplicity factor
$$  \frac{1}{(1-t^{2})(1-t^{4n})},  $$
it follows that the Betti number sequence of $P^{S^{1}}(L_{g}M;\mathbb{Z}_{2})(t)$ is unbounded. It then yields a
contradiction to Proposition \ref{Bettibound} and completes the proof. \hfill$\Box$

Next we investigate the possible B-V algebraic structures of $(\mathbb{H}_{*}(LM;\mathbb{Z}_{2}),\bullet,\bigtriangleup)$.

\begin{theorem}\label{ThmBVstr}
Let $M=\mathbb{R}P^{2n+1}$ with $n\geq1$, then  $(\mathbb{H}_{*}(LM;\mathbb{Z}_{2}),\bullet,\bigtriangleup)$ has the
following four possible B-V algebraic structures.

{\bf (A)} If $w\in\mathbb{H}_{*}(L_{e}M;\mathbb{Z}_{2})$, then we have
\be  \bigtriangleup(x)=0,\  \bigtriangleup(v)=0,\ \bigtriangleup(w)=0,  \lb{A.1}\ee
and
\be \{x,v\}=v\  {\rm or}\ v+x^{2n}vw,\ \{v,w\}=0,\ \{x,w\}=0.  \lb{A.2}\ee

{\bf (B)} If $w\in\mathbb{H}_{*}(L_{g}M;\mathbb{Z}_{2})$, then we have
\be  \bigtriangleup(x)=0,\  \bigtriangleup(v)=0,\ \bigtriangleup(w)=0,    \lb{B.1}\ee
and
\be  \{x,v\}=v,\ \{v,w\}=0,\ \{x,w\}= w\ {\rm or}\ w+x^{2n}vw^{2}.   \lb{B.2}\ee

\end{theorem}

\Proof We carry out the proof in two steps.

{\bf Step 1.} {\it Studies in Case (A).}

When $w\in\mathbb{H}_{*}(L_{e}M;\mathbb{Z}_{2})$, then $\mathbb{H}_{*}(L_{e}M;\mathbb{Z}_{2})$ is generated by
$$  \{x^{a}w^{c}\mid 0\leq a\leq 2n+1,\ c\geq 0\}.  $$
By Lemma \ref{cruciallemma}, we have
\begin{equation}\label{BV1a}
\bigtriangleup(x)=0,\ \bigtriangleup(w)=0\ \text{and}\ \bigtriangleup(xw)=0,
\end{equation}
by the B-V formulae which imply
\begin{equation}
\label{BV1b}\{x,w\}=\bigtriangleup(xw)-\bigtriangleup(x)w-x\bigtriangleup(w)=0.
\end{equation}

Here note that there is an alternative topological way to prove $\bigtriangleup(x)=0$. Let $X$ be a representative of $x$,
i.e., $[X]=x$. Since the action $\eta:S^{1}\times L_{e}M\to L_{e}M$ of (\ref{CSoperator}) is trivial on
$M\subseteq L_{e}M$ and $x\in\mathbb{H}_{-1}(M;\mathbb{Z}_{2})$, we have $\eta(S^{1}\times X)=X$ and so
$$  \bigtriangleup(x)=\eta_{*}([S^{1}]\times x)=[\eta(S^{1}\times X)]=[X]=0,  $$
where $[S^{1}]$ is the generator of $H_{1}(S^{1};\mathbb{Z}_{2})$ and the last equality is due to
$\bigtriangleup(x)\in\mathbb{H}_{0}(L_{e}M;\mathbb{Z}_{2})$.

Next we consider the fibration
\begin{equation*}\label{fiber1}
L_{e}M\to L_{e}M\times_{S^1}ES^{1}\to BS^{1},
\end{equation*}
and its Leray-Serre spectral sequence $\{\mathbb{E}^{r}, d_{r}\} =\LS(L_eM)$.

Since $d_{2}(u\otimes x^{2n+1} ) = \bigtriangleup(x^{2n+1}) = 0$ by Lemma \ref{xym1},
$x^{2n}\in\mathbb{H}_{-2n}(L_{e}M;\mathbb{Z}_{2})$ survives until $\mathbb{E}^{\infty}$. Note that the degree
of $x^{2n}$ is $-2n$, and therefore it contributes the term $t^{-2n}$ to $\mathbb{H}_{*}^{S^1}(L_{e}M;\Z_{2})$.
However, since $\mathbb{H}_{*}^{S^1}(L_{e}M;\Z_{2})$ is the desuspension of $H_{*}^{S^1}(L_{e}M;\Z_{2})$, one
should multiply it by $t^{\dim M}=t^{2n+1}$ to get a term $t$ in $P^{S^1}(L_{e}M;\mathbb{Z}_{2})(t)$ finally
(cf. Proof of Theorem 1.1 in \cite{West2007}, pp. 321-322). Here note that
the Poincar$\acute{e}$ series $P^{S^1}(L_{e}M;\Z_{2})(t)$ computed in our paper corresponds to
$H_{*}^{S^1}(L_{e}M;\Z_{2})$ instead of $\mathbb{H}_{*}^{S^1}(L_{e}M;\Z_{2})$.

However by Proposition \ref{west07}, we have
$$ P^{S^1}(LM;\mathbb{Z}_{2})(t) = \frac{1}{1-t^{2n}}\left(1+\frac{1+t}{1-t^{2}}\right)\frac{1-t^{2n+2}}{1-t^{2}},  $$
from which we know that the coefficient of $t$ in $P^{S^1}(LM;\mathbb{Z}_{2})(t)$ is $1$. Since
$$  P^{S^1}(LM;\mathbb{Z}_{2})(t) = P^{S^1}(L_{e}M;\mathbb{Z}_{2})(t)+P^{S^1}(L_{g}M;\mathbb{Z}_{2})(t),  $$
the Leray-Serre spectral sequence of the fiberation
\begin{equation*}\label{fiber2}
L_{g}M\to L_{g}M \times_{S^1}ES^{1}\to BS^{1},
\end{equation*}
must kill $x^{2n}v\in\mathbb{H}_{-2n}(L_{g}M;\mathbb{Z}_{2})$, when it passes to the third page from the second
one. As a result, we obtain $\bigtriangleup(x^{2n+1}v)=x^{2n}v$. Thus by Lemma \ref{fanfujisuan}, we obtain
\begin{equation}\label{1229a}
0 = \bigtriangleup^{2}(x^{2n+1}v) = \bigtriangleup(x^{2n}v) = x^{2n}\bigtriangleup(v),
\end{equation}
and so
\bea  x^{2n}v
&=& \bigtriangleup(x^{2n+1}v)  \nn\\
&=& \bigtriangleup(x^{2n+1})v+x^{2n+1}\bigtriangleup(v)+\{x^{2n+1},v\}  \nn\\
&=& \{x^{2n+1},v\}  \nn\\
&=& x^{2n}\{x,v\}.   \label{1229b}\eea
By (\ref{1229b}), $\{x,v\}\neq0$ and so  we can assume
$$  \{x,v\} = x^{a}v^{b}w^{c}$$
with $0\le a\le 2n+1$, $0\le b\le 1$ and $c\ge 0$, or a sum of terms like this. Then we get
$$  -a+2nc = |x^{a}v^{b}w^{c}| = |\{x,v\}| = |x|+|v|+1 = 0,  $$
i.e., $c=\frac{a}{2n}$. It then can be checked that only $(a,c)=(0,0)$ and $(a,c)=(2n,1)$ are the
possible required pairs of non-negative integers for all $n\ge 1$. Observing by the B-V formulae
$$  \{x,v\} = \bigtriangleup(x\bullet v)-\bigtriangleup(x)\bullet v-x\bullet\bigtriangleup(v),  $$
we have $\{x,v\}\in\mathbb{H}_{*}(L_{g}M;\mathbb{Z}_{2})$. Therefore, we have
$$  \{x,v\}=v,\ x^{2n}vw\ \text{or}\ v+x^{2n}vw.  $$
Again by (\ref{1229b}), we have
\be  \{x,v\} = v\  \text{or}\ v+x^{2n}vw.  \label{BV1c}\ee

Similarly if $\bigtriangleup(v)\neq 0$, we can assume that
$$  \bigtriangleup(v)=x^{a}v^{b}w^{c} $$
with $0\leq a\leq2n+1$, $0\leq b\leq1$ and $c\geq0$, or a sum of terms like this. Then we get
$$  -a+2nc = |x^{a}v^{b}w^{c}| = |\bigtriangleup(v)| = |v|+1 = 1,  $$
i.e., $c=\frac{a+1}{2n}$. It can be checked that only
$$  (a,c)=\left\{\begin{array}{ll}
              (2n-1,1),&\ \text{if}\ n\geq2;\\
              (1,1)\ \text{or}\ (3,2),&\ \text{if}\ n=1,
\end{array}\right.$$
are the possible required pairs of non-negative integers for all $n\ge 1$.

Since $\bigtriangleup:\mathbb{H}_{*}(L_{g}M;\mathbb{Z}_{2})\to\mathbb{H}_{*}(L_{g}M;\mathbb{Z}_{2})$,
by (\ref{1229a}) we have
$$\bigtriangleup(v)=\left\{\begin{array}{ll}
0\ \text{or}\ x^{2n-1}vw,&\ \text{if}\ n\geq2;\\
0\ \text{or}\ x^{3}vw^{2},&\ \text{if}\ n=1.
\end{array}\right. $$

Furthermore, we claim
\be  \bigtriangleup(v) = 0.  \lb{BV1d}\ee

In fact, we prove (\ref{BV1d}) in two cases $n\ge 2$ and $n=1$.

When $n\ge 2$, if $\bigtriangleup(v) = x^{2n-1}vw$ for some $n\geq2$, then by (\ref{BV1c}) we have
$$\begin{aligned} \bigtriangleup(xv)
&= \bigtriangleup(x)v+x\bigtriangleup(v)+\{x,v\}\\
&= x^{2n}vw+\{x,v\}\\
&= \left\{\begin{array}{ll}
       x^{2n}vw+v,&\ if\ \{x,v\}=v;\\
       v,         &\ if\ \{x,v\}=v+x^{2n}vw,\\
\end{array}\right. \end{aligned}$$
and so
$$\begin{aligned} 0 = \bigtriangleup^{2}(xv)
&=\left\{\begin{array}{ll}
  \bigtriangleup(x^{2n}vw)+\bigtriangleup(v),&\ if\ \{x,v\}=v;\\
  \bigtriangleup(v),                         &\ if\ \{x,v\}=v+x^{2n}vw,\\
\end{array}\right.\\
&= \left\{\begin{array}{ll}
  x^{2n}\bigtriangleup(vw)+x^{2n-1}vw,&\ if\ \{x,v\}=v;\\
  x^{2n-1}vw,                         &\ if\ \{x,v\}=v+x^{2n}vw,\\
\end{array}\right.\\
\end{aligned}$$
a contradiction.

When $n=1$, if $\bigtriangleup(v) = x^{3}vw^{2}$, then
$$  0 = \bigtriangleup^{2}(v) = \bigtriangleup(x^{3}vw^{2}) = \bigtriangleup(x^{3}v)w^{2} = x^{2}vw^{2},  $$
a contradiction too, and then (\ref{BV1d}) is proved.

Finally, we come to prove
\be   \{v,w\}=0.  \lb{BV1e}\ee

In fact, by arguments on topological degrees similar to the above discussion, we have
$$ \{v,w\}=\left\{\begin{array}{ll}
  0\ \text{or}\ x^{2n-1}vw^{2},                             &\ \text{if}\ n\geq2;\\
  0,\ xvw^{2},\ x^{3}vw^{3}\ \text{or}\ xvw^{2}+x^{3}vw^{3},&\ \text{if}\ n=1.
\end{array}\right.$$

If $\{v,w\}=x^{2n-1}vw^{2}$ for some $n\geq2$ or $\{v,w\}=xvw^{2}+cx^{3}vw^{3}$ with $c\in\mathbb{Z}_{2}$
when $n=1$, then  it follows from (\ref{BV1a}), (\ref{BV1c}) and (\ref{BV1d}) that
$$\begin{aligned}0=x^{2}\bigtriangleup^{2}(vw)
&=\left\{\begin{array}{ll}
x^{2}\bigtriangleup(x^{2n-1}vw^{2}),&\ \text{if}\ \{v,w\}=x^{2n-1}vw^{2};\\
x^{2}\bigtriangleup(xvw^{2}+cx^{3}vw^{3}),&\ \text{if}\ \{v,w\}=xvw^{2}+cx^{3}vw^{3},\\
\end{array}
\right.\\
&=x^{2n}\bigtriangleup(xv)w^{2}\\
&=x^{2n}\{x,v\}w^{2}\\
&=\left\{\begin{array}{ll}
x^{2n}vw^{2},&\ if\ \{x,v\}=v;\\
x^{2n}vw^{2}+x^{4n}vw^{3},&\ if\ \{x,v\}=v+x^{2n}vw,\\
\end{array}
\right.\\
&=x^{2n}vw^{2},
\end{aligned}$$
where the third identity we have used $x^{2}\bigtriangleup(cx^{3}vw^{3})=cx^{4}\bigtriangleup(xvw^{3})=0,$ a contradiction.

When $n=1$ and $\{v,w\}=x^{3}vw^{3}$, then we have
\bea 0
&=& \bigtriangleup^{2}(vw) = \bigtriangleup(x^{3}vw^{3})    \nn\\
&=& x^{2}\bigtriangleup(xvw)w^{2}  \nn\\
&=& x^{2}(\bigtriangleup(xv)w+xv\bigtriangleup(w)+\{xv,w\})w^{2}  \nn\\
&=& x^{2}(\{x,v\}w+\{xv,w\})w^{2}  \nn\\
&=& \left\{\begin{array}{ll}
     x^{2}(vw+x\{v,w\})w^{2},            &\ \text{if}\ \{x,v\}=v;\\
     x^{2}(vw+x^{2}vw^{2}+x\{v,w\})w^{2},&\ \text{if}\ \{x,v\}=v+x^{2}vw, \end{array}\right. \nn\\
&=& x^{2}vw^{3},   \label{1224g}\eea
which yields a contradiction too, and then (\ref{BV1e}) is proved.

Now by (\ref{BV1a}), (\ref{BV1b}), (\ref{BV1c}), (\ref{BV1d}) and (\ref{BV1e}), both (\ref{A.1}) and
(\ref{A.2}) are proved.

{\bf Step 2.} {\it Studies in Case (B).}

In Case (B), we have $w\in\mathbb{H}_{*}(L_{g}M;\mathbb{Z}_{2})$. Then $\mathbb{H}_{*}(L_{e}M,\mathbb{Z}_{2})$
is generated by
$$  \{x^{a}vw^{2k+1}, \ x^{a}w^{2k}\mid 1\leq a\leq 2n+1,\ k\ge 0\}.  $$
By Lemma \ref{cruciallemma}, we have
\be  \bigtriangleup(x)=0,\ \bigtriangleup(vw)=0\ \text{and}\ \bigtriangleup(xvw)=0.  \label{BV2a}\ee

From (\ref{BV2a}) we obtain
$$\begin{aligned}  0 = \bigtriangleup(xvw)
&= \bigtriangleup(x)vw+x\bigtriangleup(vw)+\{x,vw\} \\
&= \{x,vw\} \\
&= \{x,v\}w+\{x,w\}v. \end{aligned}$$
That is,
\be   \{x,v\}w = \{x,w\}v.  \label{201414a}\ee
Moreover, by the arguments almost word by word as in Case (A) we get
\be  \{x,v\}=v\ \text{or}\ v+x^{2n}w\ \text{and}\ \bigtriangleup(v)=0. \label{BV2b}\ee

By the same arguments on topological degrees as in Case (A), we have
\be  \{x,w\}=0,\ w,\ x^{2n}vw^{2}\ \text{or}\ w+x^{2n}vw^{2},  \label{201414b}\ee
and
\be   \bigtriangleup(w) = \left\{\begin{array}{ll}
  0\ \text{or}\ x^{2n-1}vw^{2},                           &\ \text{if}\ n\ge 2;\\
  0,\ xvw^{2},\ x^{3}w^{3}\ \text{or}\ xvw^{2}+x^{3}w^{3},&\ \text{if}\ n=1.  \end{array}\right. \lb{Del}\ee
From (\ref{201414a}), (\ref{BV2b}) and (\ref{201414b}) it follows
\be  \{x,v\}=v\ \text{and}\ \{x,w\}=w\ \text{or}\ w+x^{2n}vw^{2}.  \label{BV2c}\ee

We claim also
\be  \bigtriangleup(w)=0.  \label{BV2d}\ee

In fact, by (\ref{Del}) we consider two cases for $n\ge 2$ and $n=1$.

When $n\ge 2$, if $\bigtriangleup(w)=x^{2n-1}vw^{2}$, then we have
$$\begin{aligned} 0 = \bigtriangleup^{2}(w) = \bigtriangleup(x^{2n-1}vw^{2})
&= x^{2n-2}\bigtriangleup(xv)w^{2}\\
&= x^{2n-2}\{x,v\}w^{2}\\
&= \left\{\begin{array}{ll}
    x^{2n-2}vw^{2},              &\ if\ \{x,v\}=v;\\
    x^{2n-2}vw^{2}+x^{4n-2}w^{3},&\ if\ \{x,v\}=v+x^{2n}w,\\
\end{array}\right.\\
\end{aligned}$$
which yields a contradiction and proves (\ref{BV2d}) in this case.

When $n=1$ and $\bigtriangleup(w)=xvw^{2}+cx^{3}w^{3}$ with $c\in\mathbb{Z}_{2}$, an argument
similar to that for the case $n\ge 2$ yields also a contradiction.

When $n=1$ and $\bigtriangleup(w)=x^{3}w^{3}$, then we have
$$\begin{aligned} 0 = \bigtriangleup^{2}(w)
&= \bigtriangleup(x^{3}w^{3})\\
&= x^{2}\bigtriangleup(xw)w^{2}\\
&= x^{2}(\bigtriangleup(x)w+x\bigtriangleup(w)+\{x,w\})w^{2}\\
&= x^{2}\{x,w\}w^{2},\\
&= \left\{\begin{array}{ll}
      x^{2}w^{3},            &\ if\ \{x,w\}=w;\\
      x^{2}w^{3}+x^{4}vw^{4},&\ if\ \{x,w\}=w+x^{2}vw^{2}, \end{array}\right.\\
&= x^{2}w^{3},
\end{aligned}$$
which yields a contradiction too, and thus (\ref{BV2d}) is proved too.

Now it follows immediately by (\ref{BV2a}), (\ref{BV2b}) and (\ref{BV2d}) that
\be  \{v,w\} = \bigtriangleup(vw) - \bigtriangleup(v)w - v\bigtriangleup(w) = 0,  \label{BV2e}\ee
and then (\ref{B.1}) and (\ref{B.2}) are proved.

The proof of Theorem \ref{ThmBVstr} is complete. \hfill$\Box$

\setcounter{figure}{0}
\setcounter{equation}{0}
\section{Proof of Theorem \ref{Thm1.1}}%Section 5

Our idea of the proof of Theorem \ref{Thm1.1} is simple. We compute first the Poincar\'e series associated to the
third pages of the Leray-Serre spectral sequences $\LS(L_{e}M)$ and $\LS(L_{g}M)$ for each of the four possible
B-V algebraic structures obtained in Theorem \ref{ThmBVstr} respectively. Then Theorem \ref{Thm1.1} follows by
comparing their sums with the result obtained by Westerland in \cite{West2007}, i.e., Proposition \ref{west07}.

{\bf Proof of Theorem \ref{Thm1.1}:}

We carry out the proof in two Cases (A) and (B) according to Theorem \ref{ThmBVstr}.

{\bf Step 1.} {\it Studies in Case (A) of Theorem \ref{ThmBVstr}.}

In this case, we have $w\in\mathbb{H}_{*}(L_{e}M;\mathbb{Z}_{2})$, and by (\ref{A.1}) and (\ref{A.2}) of
Theorem \ref{ThmBVstr}, we have
\bea
&& \bigtriangleup(x)=0,\ \bigtriangleup(v)=0\ \text{and}\ \bigtriangleup(w)=0,  \nn\\
&& \{x,v\}=v\ \text{or}\ v+x^{2n}vw,\ \{x,w\}=0\ \text{and}\ \{v,w\}=0.  \nn\eea
We continue in three sub-steps.

{\bf (i)} {\it Now we consider the fibration
$$   L_{e}M\to L_{e}M\times_{S^1}ES^{1}\to BS^{1},  $$
and its Leray-Serre spectral sequence $\LS(L_{e}M)$.}

Since $\bigtriangleup \equiv 0$ on $\mathbb{H}_{*}(L_{e}M;\mathbb{Z}_{2})$ by Lemma \ref{cruciallemma}, we have
$d_{2}\equiv 0$ on $\mathbb{E}^{2}$ by Lemma \ref{xym1}. Therefore $\mathbb{E}^{3}$ is the same as
$\mathbb{E}^{2}$ and so the Poincar\'e series $P_{III}^{S^1}(L_{e}M;\mathbb{Z}_{2})(t)$ of the third page
equals to that of the second page. Observing that
$$  \mathbb{E}^{2}=\mathbb{H}_{*}(L_{e}M)\otimes H_{*}(BS^{1})\cong\mathbb{Z}_{2}[x,w,u]/(x^{2n+2}),  $$
with $|x|=-1$, $|w|=2n$ and $|u|=2$ given in Lemma \ref{xym1}, we get
\bea  P_{III}^{S^1}(L_{e}M;\mathbb{Z}_{2})(t)
&=& P_{II}^{S^1}(L_{e}M;\mathbb{Z}_{2})(t)   \nn\\
&=& t^{2n+1}\sum_{a=0}^{2n+1}\sum_{b=0}^{+\infty}\sum_{c=0}^{+\infty}t^{-a}t^{2b}t^{2nc}   \nn\\
&=& t^{2n+1} \left(\frac{1-t^{-(2n+2)}}{1-t^{-1}}\right)\frac{1}{(1-t^{2})(1-t^{2n})}  \nn\\
&=& \frac{1}{1-t^{2n}}\left(\frac{1-t^{2n+2}}{1-t^{2}}\right) \frac{1+t}{1-t^{2}}.  \label{le3}\eea

{\bf (ii)} {\it We consider the fibration
$$  L_{g}M\to L_{g}M \times_{S^1}ES^{1}\to BS^{1},  $$
and its Leray-Serre spectral sequence $\LS(L_{g}M)$.}

By direct computations on $\mathbb{H}_{*}(L_{g}M;\mathbb{Z}_{2})$ we obtain
\bea  \bigtriangleup(x^{2l}vw^{c})
&=& x^{2l}\bigtriangleup(vw^{c})   \nn\\
&=& x^{2l}(\bigtriangleup(v)w^{c}+v\bigtriangleup(w^{c})+\{v,w^{c}\})   \nn\\
&=& 0,  \nn\eea
and
\bea \bigtriangleup(x^{2l+1}vw^{c})
&=& x^{2l}\bigtriangleup(xvw^{c})   \nn\\
&=& x^{2l}(\bigtriangleup(xv)w^{c} + xv\bigtriangleup(w^{c}) + \{xv,w^{c}\}) \nn\\
&=& x^{2l}\bigtriangleup(xv)w^{c} \nn\\
&=& x^{2l}\{x,v\}w^{c}  \nn\\
&=& \left\{\begin{array}{ll}
    x^{2l}vw^{c},                  &\ \text{if}\ \{x,v\}=v;  \\
    x^{2l}vw^{c}+x^{2l+2n}vw^{c+1},&\ \text{if}\ \{x,v\}=v+x^{2n}vw,  \\
    \end{array}\right.  \nn\\
&=& \left\{\begin{array}{ll}
    x^{2l}vw^{c},         &\ \text{if}\ \{x,v\}=v;\\
    x^{2l}vw^{c},         &\ \text{if}\ \{x,v\}=v+x^{2n}vw\ \text{and}\ l\geq1,\\
    vw^{c}+x^{2n}vw^{c+1},&\ \text{if}\ \{x,v\}=v+x^{2n}vw\ \text{and}\ l=0,
    \end{array}\right.  \lb{5.2}\eea
for $0\le l\le n$ and $0\le c< +\infty$.

For the case $\{x,v\}=v+x^{2n}vw$, by (\ref{5.2}) we have also
\bea \bigtriangleup(xvw^{c}+x^{2n+1}vw^{c+1})
&=& \bigtriangleup(xvw^{c})+\bigtriangleup(x^{2n+1}vw^{c+1})   \nn\\
&=& (vw^{c}+x^{2n}vw^{c+1})+x^{2n}vw^{c+1}\\&=&vw^{c}.   \nn\eea

So no matter $\{x,v\}=v$ or $v+x^{2n}vw$, only the elements generated by
$$  \{x^{2l+1}vw^{c}\mid 0\le l\le n,\ c\ge 0\}  $$
survive when $\LS(L_{g}M)$ passes to the third page from the second one, while the other elements generated by
$$  \{x^{2l}vw^{c}\mid 0\le l\le n,\ c\ge 0\}\cup\{u^{p}\otimes x^{a}vw^{c} \mid p\ge 1, 0\le a\le 2n+1,\ c\ge 0\}, $$
are killed since they are either in the image or not in the kernel of the second differential. As a result, we obtain
\bea   P_{III}^{S^1}(L_{g}M;\mathbb{Z}_{2})(t)
&=& t^{2n+1}\sum_{l=0}^{n}\sum_{c=0}^{+\infty}t^{-(2l+1)}t^{2nc}    \nn\\
&=& t^{2n} \left(\frac{1-t^{-(2n+2)}}{1-t^{-2}}\right) \frac{1}{1-t^{2n}}    \nn\\
&=& \frac{1}{1-t^{2n}}\left(\frac{1-t^{2n+2}}{1-t^{2}}\right).  \label{lg3}\eea

{\bf (iii)} {\it Conclusion on $P^{S^1}(L_{g}M;\mathbb{Z}_{2})(t)$.}

Comparing the sum of the two Poincar\'e series (\ref{le3}) and (\ref{lg3}) with $P^{S^1}(LM;\mathbb{Z}_{2})(t)$ given
by Proposition \ref{west07}, we have
$$ P_{III}^{S^1}(L_{e}M;\mathbb{Z}_{2})(t)+P_{III}^{S^1}(L_{g}M;\mathbb{Z}_{2})(t)=P^{S^1}(LM;\mathbb{Z}_{2})(t), $$
which implies that $\LS(L_{e}M)$ and $\LS(L_{g}M)$ collapse at the second and third pages respectively. Thus finally
we obtain
$$ P^{S^1}(L_{g}M;\mathbb{Z}_{2})(t) = P_{III}^{S^1}(L_{g}M;\mathbb{Z}_{2})(t)
     = \frac{1}{1-t^{2n}}\left(\frac{1-t^{2n+2}}{1-t^{2}}\right),  $$
which is precisely (\ref{Clm1.1}) of Theorem \ref{Thm1.1} in Case (A).

{\bf Step 2.} {\it Studies in Case (B) of Theorem \ref{ThmBVstr}.}

In this case, $w\in\mathbb{H}_{*}(L_{g}M;\mathbb{Z}_{2})$. Then by (\ref{B.1}) and (\ref{B.2}) in Theorem \ref{ThmBVstr}
we have
\bea
&& \bigtriangleup(x)=0,\ \bigtriangleup(v)=0\ \text{and}\ \bigtriangleup(w)=0,  \nn\\
&& \{x,v\}=v,\ \{x,w\}=w\ \text{or}\ w+x^{2n}vw^{2}\ \text{and}\ \{v,w\}=0.  \nn\eea

By similar arguments as in Step 1, we get
\bea  P_{III}^{S^1}(L_{e}M;\mathbb{Z}_{2})(t)
&=& P_{II}^{S^1}(L_{e}M;\mathbb{Z}_{2})(t)  \nn\\
&=& \frac{1}{1-t^{2n}} \left(\frac{1-t^{2n+2}}{1-t^{2}}\right) \frac{1+t}{1-t^{2}}. \nn\eea

On the other hand, by direct computations on $\mathbb{H}_{*}(L_{g}M;\mathbb{Z}_{2})$ we obtain

(i) if $\{x,w\}=w$, then
$$ \bigtriangleup(x^{a}v^{b}w^{c})
=\left\{\begin{array}{ll}
    0,                &\ \text{if}\ a\ \text{is\ even};\\
    x^{a-1}v^{b}w^{c},&\ \text{if}\ a\ \text{is\ odd},\\
\end{array}\right.$$
and

(ii) if $\{x,w\}=w+x^{2n}vw^{2}$, then
$$  \bigtriangleup(x^{a}v^{b}w^{c})
= \left\{\begin{array}{ll}
    0,                              &\ \text{if}\ a\ \text{is\ even};\\
    x^{a-1}vw^{c},                  &\ \text{if}\ a\ \text{is\ odd\ and}\ c\ \text{is\ even};\\
    x^{a-1}w^{c}+x^{a-1+2n}vw^{c+1},&\ \text{if}\ a\ \text{is\ odd\ and}\ c\ \text{is\ odd},\\
\end{array}\right.$$
where
$$ b=\left\{\begin{array}{ll}
    1,&\ \text{when}\ c\ \text{is even};\\
    0,&\ \text{when}\ c\ \text{is odd}.
\end{array}\right.$$

For the latter case, we also have
\bea  \bigtriangleup(xw^{c}+x^{2n+1}vw^{c+1})
&=& \bigtriangleup(xw^{c})+\bigtriangleup(x^{2n+1}vw^{c+1})  \nn\\
&=& (w^{c}+x^{2n}vw^{c+1})+x^{2n}vw^{c+1}    \nn\\
&=& w^{c},   \nn\eea
for odd $c$.

Thus no matter $\{x,w\}=w$\ or $w+x^{2n}vw$, only the elements generated by
$$  \{x^{2l+1}w^{2k+1}\mid 0\le l\le n, \;k\ge 0\},  $$
and
$$  \{x^{2l+1}vw^{2k}\mid 0\le l\le n, k\ge 0\},  $$
survive when $\LS(L_{g}M)$ passes to the third page from the second one. As a result, we obtain
\bea   P_{III}^{S^1}(L_{g}M;\mathbb{Z}_{2})(t)
&=& t^{2n+1}\sum_{l=0}^{n}\sum_{k=0}^{+\infty}t^{-2l-1}t^{4kn+2n}
           + t^{2n+1}\sum_{l=0}^{n}\sum_{k=0}^{+\infty}t^{-2l-1}t^{4kn}  \nn\\
&=& t^{2n+1}\sum_{l=0}^{n}\sum_{k=0}^{+\infty}t^{-2l-1}t^{2kn}  \nn\\
&=& \frac{1}{1-t^{2n}} \left(\frac{1-t^{2n+2}}{1-t^{2}}\right).  \nn\eea

The rest proof is then word by word as that in Step 1.

Finally using (\ref{Clm1.1}) of $P^{S^{1}}(\Lambda_{g}M;\mathbb{Z}_{2})(t)$ to find the coefficients
$\bar{\beta}_{i}$s, by direct computations we obtain $\bar{B}_g(M)= (n+1)/(2n)$, i.e., (\ref{aB.2})
holds, and complete the proof of the Theorem \ref{Thm1.1}. \hfill$\Box$

\begin{remark}Theorem \ref{Thm1.1} is not a trivial conclusion of Lemma \ref{cruciallemma} (or Theorem
\ref{ThmBVstr}). In fact, if the sum of $P_{III}^{S^1}(L_{e}M;\mathbb{Z}_{2})(t)$ and
$P_{III}^{S^1}(L_{g}M;\mathbb{Z}_{2})(t)$ was greater than $P^{S^1}(LM;\mathbb{Z}_{2})(t)$ under one
of the four possible B-V algebraic structures, one can not claim  by Proposition \ref{Bettibound} that
the same conclusion of Lemma \ref{cruciallemma} holds for the higher even differentials of $\LS(L_{e}M$).
 The reason is that we do not know whether the higher even differentials of $\LS(L_{e}M$)
have the ``homogeneous'' property possessed by the second one (Lemma \ref{xym1}), i.e., either it keeps
the associated page stable or kills a ``large'' series, whose Betti number sequence is
unbounded (cf. the proof of Lemma \ref{cruciallemma}).

But it is then difficult to know at which pages $\LS(L_{e}M)$ and $\LS(L_{g}M)$ collapse because the higher
even differentials of the spectral sequences go mysteriously, provided that  the above mentioned phenomenon
happened. As a result, we could not obtain $P^{S^1}(L_{g}M;\Z_{2})(t)$ any more.
\end{remark}

\setcounter{figure}{0}
\setcounter{equation}{0}
\section{Proof of Theorem \ref{Thm1.2}}%Section 6

In this section, we apply Theorem \ref{Thm1.1} to obtain the resonance identity of non-contractible prime
homologically visible prime closed geodesics on a Finsler $M=(\mathbb{R}P^{2n+1},F)$ claimed in Theorem
\ref{Thm1.2}, provided the number of all the distinct prime closed geodesics on $M$ is finite.

{\bf The proof of Theorem \ref{Thm1.2}}. Recall that we denote the homologically visible prime closed
geodesics by $\CG_{\hv}(M)=\{c_1, \ldots, c_r\}$ for some integer $r>0$ when the number of distinct prime
closed geodesics on $M$ is finite. Note also that by Lemma \ref{Lm3.4} we have $\hat{i}(c_j)>0$ for all
$1\le j\le r$.

Let $w_h = M_h(\Lm_gM)$ defined by (\ref{M.1}). The Morse series of $\Lm_gM$ is defined by
\be  M(t) = \sum_{h=0}^{+\infty}w_ht^h.  \label{wh}\ee
Note that $\{w_h\}$ is a bounded sequence by the second inequality of (\ref{M.2}). We now use the method
in the proof of Theorem 5.4 of \cite{LW2007} to estimate
$$   M^{q}(-1) = \sum_{h=0}^{q}w_h(-1)^h.   $$

By (\ref{wh}) and (\ref{CGap2}) we obtain
\bea M^{q}(-1)
&=& \sum_{h=0}^{q}w_{h}(-1)^{h}  \nn\\
&=& \sum_{j=1}^{r}\sum_{m=1}^{n_j/2}\sum_{l=0}^{4n}\sum_{h=0}^{q}(-1)^{h}k_{l}(c_{j}^{2m-1})
            \;{}^{\#}\left\{s\in\mathbb{N}\cup\{0\}\mid h-i(c_{j}^{2m-1+sn_j})=l\right\}  \nn\\
&=& \sum_{j=1}^{r}\sum_{m=1}^{n_j/2}\sum_{l=0}^{4n}(-1)^{l+i(c_{j})}k_{l}(c_{j}^{2m-1})
            \;{}^{\#}\left\{s\in\mathbb{N}\cup\{0\}\mid l+i(c_{j}^{2m-1+sn_j})\le q\right\}.  \nn\eea
On the one hand, we have
\bea
&&{}^{\#}\left\{s\in\mathbb{N}\cup\{0\}\mid l+i(c_{j}^{2m-1+sn_j})\le q\right\} \nn\\
&&\qquad = \;{}^{\#}\left\{s\in\mathbb{N}\cup\{0\}\mid l+i(c_{j}^{2m-1+sn_j})\le q,\;
           |i(c_{j}^{2m-1+sn_j})-(2m-1+sn_j)\hat{i}(c_{j})|\le 2n\right\}  \nn\\
&&\qquad \le \;{}^{\#}\left\{s\in\mathbb{N}\cup\{0\}\mid 0\le (2m-1+sn_{j})\hat{i}(c_{j})\le q-l+2n  \right\}  \nn\\
&&\qquad = \;{}^{\#}\left\{s\in\mathbb{N}\cup\{0\}\mid 0 \leq s
       \le \frac{q-l+2n-(2m-1)\hat{i}(c_j)}{n_j\hat{i}(c_j)}\right\} \nn\\
&&\qquad \le \frac{q-l+2n}{n_j\hat{i}(c_j)} + 1. \nn\eea
On the other hand, we have
\bea
&&{}^{\#}\left\{s\in\mathbb{N}\cup\{0\}\mid l+i(c_{j}^{2m-1+sn_j})\le q \right\} \nn\\
&&\qquad = \;{}^{\#}\left\{s\in\mathbb{N}\cup\{0\}\mid l+i(c_{j}^{2m-1+sn_j})\le q,\;
           |i(c_{j}^{2m-1+sn_j})-(2m-1+sn_{j})\hat{i}(c_{j})|\le 2n\right\}  \nn\\
&&\qquad \ge \;{}^{\#}\left\{s\in\mathbb{N}\cup\{0\}\mid i(c_{j}^{2m-1+sn_j})
       \le (2m-1+sn_{j})\hat{i}(c_{j})+2n\le q-l \right\} \nn\\
&&\qquad \ge \;{}^{\#}\left\{s\in\mathbb{N}\cup\{0\}\mid 0 \le s
       \le \frac{q-l-2n-(2m-1)\hat{i}(c_{j})}{n_j\hat{i}(c_{j})} \right\}  \nn\\
&&\qquad \ge \frac{q-l-2n}{n_j\hat{i}(c_{j})} - 1.  \nn\eea
Thus we obtain
$$  \lim_{q\to+\infty}\frac{1}{q}M^{q}(-1)
  = \sum_{j=1}^{r}\sum_{m=1}^{n_j/2}\sum_{l=0}^{4n}(-1)^{l+i(c_{j})}k_{l}(c_{j}^{2m-1})\frac{1}{n_{j}\hat{i}(c_{j})}
              = \sum_{j=1}^{r}\frac{\hat{\chi}(c_j)}{\hat{i}(c_j)}.  $$
Since $w_{h}$ is bounded, we then obtain
$$  \lim_{q\to+\infty}\frac{1}{q}M^{q}(-1) = \lim_{q\to+\infty}\frac{1}{q}P^{S^{1},q}(\Lambda_{g}M;\mathbb{Z}_{2})(-1)
          = \lim_{q\to+\infty}\frac{1}{q}\sum_{k=0}^{q}(-1)^k\bar{\beta}_k(\Lm_gM) = \bar{B}(\Lm_gM),  $$
where $P^{S^{1},q}(\Lambda_{g}M;\mathbb{Z}_{2})(t)$ is the truncated polynomial of
$P^{S^{1}}(\Lambda_{g}M;\mathbb{Z}_{2})(t)$ with terms of degree less than or equal to $q$.
Thus we get
$$  \sum_{j=1}^{r}\frac{\hat{\chi}(c_j)}{\hat{i}(c_j)} = \frac{n+1}{2n}, $$
which proves Theorem \ref{Thm1.2}. $\hfill\Box$

\smallskip

{\noindent\bf Acknowledgements} The authors would like to sincerely thank Professor Craig Westerland for
many helpful discussions on his works. The authors would also like to thank sincerely the referee for his/her
careful reading and valuable comments on the paper.

\end{document}